\documentclass{elsart}
\usepackage{amsmath,amssymb,epic,graphicx}
\makeatletter

\@addtoreset{equation}{section}
\makeatother

\newcommand{\id}{{1\!\!1}}
\newcommand{\Kk}{\mathbb{K}}
\newcommand{\Ff}{\mathbb{F}}
\newcommand{\Ss}{\mathbb{S}}
\newcommand{\ints}{\mathbb{Z}}
\newcommand{\reals}{\mathbb{R}}
\newcommand{\rats}{\mathbb{Q}}
\newcommand{\cx}{\mathbb{C}}
\newcommand{\quat}{\mathbb{H}}
\newcommand{\oct}{\mathbb{O}}
\newcommand{\capX}{{\tt X}}
\newcommand{\bX}{\overline{\capX}}
\newcommand{\capY}{{\tt Y}}

\newcommand{\bM}{\overline{M}}

\renewcommand{\a}{\alpha}
\newcommand{\ve}{\varepsilon}
\newcommand{\bve}{\bar{\varepsilon}}
\newcommand{\bth}{\overline{\theta}}
\newcommand{\E}{\overline{\theta}}
\newcommand{\si}{\sigma}
\newcommand{\rmD}{{\rm{Det}}}

\newcommand{\be}{\begin{eqnarray}}
\newcommand{\ee}{\end{eqnarray}}
\newcommand{\nn}{\nonumber}

\newcommand{\ca}{A}
\newcommand{\cb}{B}
\newcommand{\cc}{C}
\newcommand{\cd}{D}
\newcommand{\ce}{E}
\newcommand{\cZ}{{\mathcal{Z}}}
\newcommand{\Units}{{\mathcal E}}

\newcommand{\cN}{{\mathcal N}}
\newcommand{\cF}{{\mathcal F}}
\newcommand{\tS}{{\tilde{S}}}
\newcommand{\lp}{\left[}
\newcommand{\rp}{\right]}
\newcommand{\Lq}{\mathtt{L}}
\newcommand{\Hq}{\mathtt{H}}
\newcommand{\Iq}{\mathtt{I}}
\newcommand{\Cq}{\mathtt{C}}
\newcommand{\Ztwo}{{\mathbb{Z}}_2}
\newcommand{\cE}{\mathtt{E}}
\newcommand{\cG}{\mathtt{G}}
\newcommand{\Oo}{\mathtt{O}}
\newcommand{\cI}{\mathtt{I}}
\newcommand{\cO}{{\mathcal{O}}}
\newcommand{\cR}{{\mathtt{R}}}
\newcommand{\cT}{{\mathcal{T}}}

\newcommand{\mfg}{\mathfrak{g}}
\newcommand{\Sg}{{\mathfrak S}}
\newcommand{\Ag}{{\mathfrak A}}
\newcommand{\PGL}{{PSL_2^{(0)}(\Hq)}}
\newcommand{\tPGL}{{\widetilde{PSL}_2^{(0)}(\Hq)}}
\newcommand{\bd}{\mbox{\boldmath$.$\unboldmath}}

\newtheorem{lemma}{Lemma}

\begin{document}

{\flushright ULB-TH/08-10\\
AEI-2008-013\\[20mm]}

\begin{center}
{\LARGE \sc Hyperbolic Weyl groups and \\[4mm]
 the four normed division algebras}
\end{center}
\vskip10mm
\begin{center}
Alex J. Feingold\footnotemark[1], Axel Kleinschmidt\footnotemark[2] and
Hermann Nicolai\footnotemark[3]\\[5mm]
\footnotemark[1]{\sl  Department of Mathematical Sciences, The State
University of New York\\ 
 Binghamton, New York 13902--6000, U.S.A.}\\[3mm]
\footnotemark[2]{\sl Physique Th\'eorique et Math\'ematique,\\
Universit\'e Libre de Bruxelles \&{} International Solvay Institutes,\\
Boulevard du Triomphe, ULB -- CP 231, B-1050 Bruxelles, Belgium}\\[3mm]
\footnotemark[3]{\sl  Max-Planck-Institut f\"ur Gravitationsphysik,
 Albert-Einstein-Institut \\
 Am M\"uhlenberg 1, D-14476 Potsdam, Germany} \\[7mm]

\begin{tabular}{p{11cm}}
\hspace{5mm}{\footnotesize {\bf Abstract:} We study the Weyl groups of 
hyperbolic Kac--Moody algebras of `over-extended' type and ranks 3, 4, 6
and 10, which are intimately linked with the four normed division algebras 
$\Kk=\reals,\cx,\quat,\oct$, respectively. A crucial role is played by 
integral lattices of the division algebras and associated discrete matrix 
groups. Our findings can be summarized by saying that the even subgroups, 
$W^+$, of the Kac--Moody Weyl groups, $W$, are isomorphic to generalized 
modular groups over $\Kk$ for the simply laced algebras, and to certain 
finite extensions thereof for the non-simply laced algebras.
This hints at an extended theory of modular forms and functions.}

\end{tabular}\\[5mm]
\end{center}

\vskip 20pt

\begin{center}
``The mathematical universe is inhabited not only by important species
but also by interesting individuals.'' -- C. L. Siegel 
\end{center} 
\vskip 10pt 

\thispagestyle{empty}

\newpage
\setcounter{page}{1}

\begin{section}{Introduction}

In \cite{FeFr83} Feingold and Frenkel gained significant new insight 
into the structure of a particularly interesting rank $3$ hyperbolic 
Kac--Moody algebra which they called $\cF$, along with some connections 
to the theory of Siegel modular forms of genus 2. The first vital step 
in their work was the discovery that the Weyl group of that hyperbolic 
algebra is $W(\cF) \cong PGL_2(\ints)$, the projective group of $(2\times 2)$ 
integral matrices with determinant $\pm 1$, isomorphic to the hyperbolic 
triangle group $T(2,3,\infty)$. They showed that the root system of that 
algebra could be realized as the set of $(2\times 2)$ symmetric integral 
matrices $\capX$ with $\det(\capX) \geq -1$, and that the action of
$M\in W(\cF)$  on $\capX$ is given by $M\capX M^T$. In notation used
commonly by physicists today,  
the hyperbolic algebra studied in \cite{FeFr83} is designated as $A_1^{++}$ 
since it is obtained from the finite-dimensional Lie algebra $sl_2$ of 
type $A_1$ by a process of double extension. The first step of the extension 
gives the affine algebra $A_1^{(1)}\equiv A_1^+$, and the second step, often
referred to as over-extension, gives the hyperbolic algebra. 
In the realization of the root lattice via symmetric matrices $\capX$,
the real roots consist of the integral points $\capX$ with
$\det(\capX) = -1$ on a single-sheeted hyperboloid, and the  
imaginary roots consist of the integral points $\capX$ on the light-cone 
$\det(\capX) = 0$ and on the two-sheeted hyperboloids $\det(\capX) > 0$. 
In \cite{FeFr83} it was also mentioned that these results could be extended  
to two other (dual) rank $4$ hyperbolic algebras whose Weyl groups were 
both the Klein--Fricke group $\Psi_1^*$ containing as an index $4$ subgroup 
the Picard group $\Psi_1 = PSL_2(\mathbb{Z}(i))$ whose entries are from 
the Gaussian integers. The results expected to hold for these rank 
$4$ hyperbolic algebras included connections to the theory of Hermitian 
modular forms, but that line of research was not pursued in later work.
A paper by Kac, Moody and Wakimoto \cite{KMW} generalized the structural 
results of \cite{FeFr83} to the hyperbolic algebra $E_{10} = E_8^{++}$, 
but until now there has not been any new insight into the structure of 
the Weyl groups of hyperbolic algebras which are usually just given 
as Coxeter groups defined by generators and relations.

In the present work, we take up this line of development again, and
present a coherent picture for many higher rank hyperbolic Kac--Moody
algebras which is based on the relation to generalized modular groups
associated with lattices and subrings of the four normed division 
algebras. More specifically, we shall show that the Weyl groups of 
all hyperbolic algebras of ranks 4, 6 and 10 which can be obtained 
by the process of double extension described above, admit realizations 
in terms of generalized modular groups over the complex numbers $\cx$, 
the quaternions $\quat$, and the octonions $\oct$, respectively. 
We are encouraged to find that these Weyl groups are amenable to 
explicit matrix descriptions, but understand that the hyperbolic algebras
themselves have still eluded any effective characterization. 
For $\Kk=\cx$, the present work on these hyperbolic Weyl groups is thus a 
very natural extension of \cite{FeFr83}, showing in particular that 
the rank $4$ hyperbolics $A_2^{++}, C_2^{++}$ and $G_2^{++}$ are 
naturally connected with certain subrings in the normed division 
algebra $\cx$. Analogous results are obtained for all the `over-extended' 
rank $6$ hyperbolics $A_4^{++}, B_4^{++}, C_4^{++}, D_4^{++}$ and 
$F_4^{++}$, whose (even) Weyl groups can be described in terms of
quaternionic modular groups. Finally, $\Kk=\oct$, the largest and 
non-associative division algebra of octonions, is associated with  
rank $10$ hyperbolics, and in particular, the maximally extended
hyperbolic Kac--Moody algebra $E_{10}$. The other two hyperbolic over-extended
algebras $B_8^{++}$ and $D_8^{++}$ can also be described using octonions. 
In this paper we present the 
rich structure which we found in the  complex and quaternionic cases, 
as well as partial (and intriguing) results for the octonionic case. 
As explained in more detail in Section~\ref{weylsec}, a new feature for 
the division algebras beyond $\reals$ is that simple reflections involve 
complex conjugation of all entries of $\capX$. For that reason only the 
{\em even} part of a given Weyl group will act by matrix conjugation 
of $\capX$ and is therefore the main focus of our study. Throughout 
the paper we will denote by $W^+$ the even part of a Weyl group $W$.
In the case of $\cF$ one has $W^+(\cF)\cong PSL_2(\ints)$, the modular 
group, and for the other division algebras we discover a number of 
apparently new modular groups. An announcement of our work appears in
\cite{FeKlNi09}, based on a talk presented by A.J.F. at the conference 
on ``Vertex Operator Algebras and Related Areas" in honor of Geoffrey
Mason, July 2008, at Illinois State University. 

In table~\ref{integers} we summarize our findings for the (even) Weyl 
groups of the finite and hyperbolic Kac--Moody algebras studied in 
relation to the various division algebras. 
By $\mfg^{++}$ we mean 
the over-extension of the finite-dimensional Lie algebra $\mfg$ to
a hyperbolic Kac--Moody algebra, where the first extension $\mfg^{+}$
is an untwisted affine algebra. In Appendix A we discuss
other cases where the first extension is twisted affine, so a different
notation is required. In the table, we use the standard group theory 
notation $C=A\,\bd\,B$ to mean a group $C$ which contains group $A$ as a 
normal subgroup, with quotient group $C/A$ isomorphic to $B$. 
Such a group $C$ is called an extension of $A$ by $B$. If $B$ has order $|B|$, 
$A$ is said to be of index $|B|$ in $A\,\bd\,B$. It can happen that the
extension  is a semi-direct product, so that $B$ is a subgroup of $C$ which
acts on $A$ via  conjugation as automorphisms, and in this case the product is
denoted by  $A\,\rtimes\,B$. By $\Sg_n$ we denote the symmetric group on $n$
letters.

\begin{table}
\begin{center}
 \begin{tabular}{|c|c|c|c|}
	\hline
	\;\;\;$\Kk$\;\;\; & \;Root system\; $\mfg$ 
       &  $W(\mfg)$ &\;\;\;\qquad $W^+(\mfg^{++})$\;\;\;\qquad\\
	\hline
	$\reals$ & $A_1$ 
           & $2\equiv\Ztwo$ & $PSL_2(\ints)$\\
	$\cx$ & $A_2$ 
  &  $\ints_3 \rtimes 2$ & $PSL_2(\cE)$\\
	$\cx$ & $B_2\equiv C_2$ 
       & $\ints_4\rtimes 2$ & 
       $PSL_2(\cG)\rtimes 2$\\
	$\cx$ & $G_2$ 
       & $\ints_6 \rtimes 2$ & $PSL_2(\cE)\rtimes2$ \\
        $\quat$ & $A_4$ 
        & $\Sg_5$ &  $ PSL_2 ^{(0)}(\Iq)$\\
        $\quat$ & $B_4$ 
        & $2^4 \rtimes \Sg_4$ &  $\PGL\rtimes 2$\\
        $\quat$ & $C_4$ 
        & $2^4 \rtimes \Sg_4$ & $\tPGL \rtimes 2$\\
	$\quat$ & $D_4$ 
       & $2^3\rtimes\Sg_4$ & $\PGL$ \\
        $\quat$ & $F_4$ 
        & $2^5\rtimes (\Sg_3\times \Sg_3)$ &
           $PSL_2(\Hq)\rtimes 2$\\
        $\oct$ & $D_8$ 
       & $2^7 \rtimes \Sg_8$ & $PSL_2^{(0)}(\Oo)$\\
        $\oct$ & $B_8$ 
       & $2^8\rtimes \Sg_8$  & $PSL_2^{(0)}(\Oo)\rtimes 2$\\
	$\oct$ & $E_8$ 
       & $2\,\bd\,O_8^+(2)\,\bd\,2$ & $PSL_2(\Oo)$\\
	\hline
	\end{tabular}
	\caption{ \label{integers}\em Root systems, number systems and Weyl
	groups within the normed division algebras. The groups in the
       right column are defined in detail in the relevant sections.} 
\end{center}
\end{table}

The various rings appearing in the table are as follows. 
For $\Kk=\cx$ these are the Gaussian integers
$\cG \equiv \ints(i)$ and the Eisenstein integers $\cE$; for the
quaternions $\quat$ we use the terminology of \cite{CoSm03}, referring 
to the maximal order of quaternions
having all coefficients in $\ints$ or all in $\ints\! +\! \frac{1}{2}$ as 
{\em Hurwitz}  integers  $\Hq$. To understand $W^+(A_4^{++})$ we will need
the icosian quaternions $\Iq$ \cite{CoSl88,MoPa93}. 
The octonionic integers are called the {\em octavians} $\Oo$ \cite{CoSm03}.
The generalized modular groups that we find are all discrete as matrix groups.

Our results are complete for the cases $\Kk=\cx$ and $\Kk=\quat$.
For all of these we can reformulate the even Weyl groups as new kinds 
of modular groups, most of which have not yet appeared in the literature 
so far. We note that A.~Krieg \cite{Kr85,Krieg1} has developed a 
theory of modular forms over the quaternions, however the modular 
groups defined there do not appear to coincide with the ones found here.
Nevertheless, we will make extensive use of some results of \cite{Kr85} 
in our analysis, in particular Theorem~2.2 on p.~16. 

The most interesting Weyl group, however, is the one of the maximally 
extended hyperbolic Kac--Moody algebra $E_{10} \equiv E_8^{++}$, where 
our results are still incomplete, and which provided the chief motivation 
for the present work. In this case, the even Weyl group is given by 
\be
W^+(E_{10}) \cong PSL_2(\Oo)\,.
\ee 
Because the $(2\times 2)$ matrices over the octonions obviously do not form
a group, the `modular group' over the octavians $\Oo$ appearing on the
r.h.s. of this equation can so far only be defined recursively, by nested 
sequences of matrix conjugations with the generating (octonionic)
matrices introduced in Section~\ref{weylsec} of this paper. This is in
analogy with  
the generators and relations description of the $E_{10}$ Weyl group as products 
of fundamental reflections (and also in analogy with the description of 
the continuous Lorentz group $SO(1,9;\reals)$ via octonionic $(2\times 2)$
matrices  in \cite{Su84,Ma93}). It remains an outstanding problem to find a
more  manageable realization of this modular group directly in terms of 
$(2\times 2)$ matrices with octavian entries and the $G_2(2)$ automorphism group
of the octavian integers. In the final section we present some results 
in this direction which we believe are new, and which may be of use in 
future investigations. The case $E_{10}$ is also the most important 
because it has recently been shown that all simply laced hyperbolic 
Kac--Moody algebras can be embedded into $E_{10}$ \cite{Vi08}. Because 
the associated Weyl groups are then subgroups of $W(E_{10})$, all 
these Weyl subgroups should admit an octonionic realization. Conversely, 
the structure and explicit realizations of modular groups found 
here for the quaternionic case should help in understanding the
group $PSL_2(\Oo)$, because all rank 6 algebras occur as subalgebras
inside $E_{10}$ (and $A_4^{++}$ and $D_4^{++}$ as regular subalgebras,
in particular). The remarks above also apply to the Weyl groups of
$B_8^{++}$ and $D_8^{++}$.

Ref.~\cite{FeFr83} also highlighted a possible link with the theory of
{\em Siegel modular forms}. In the present context, one would thus
start from the Jordan algebra $H_2(\Kk_\cx)$ over the complexified
division algebra $\Kk_\cx$, and consider a generalized Siegel upper
half-plane for $\cZ\in H_2(\Kk_\cx)$. This can be done not only for 
the quaternions \cite{Krieg1}, but also for $\Kk=\oct$, where the 
outlines of a corresponding theory of `Siegel modular forms' have been 
described and developed in \cite{Krieg2}. However, it remains to be seen 
what role $PSL_2(\Oo)$ as defined here has to play in this context, 
and whether results in this direction can substantially advance our 
understanding of the hyperbolic algebra $E_{10}$.

We believe that our results are also interesting in the light of
recent developments in the study of black hole microstates for
particular classes of black holes. It was found already in \cite{DVV}
that the degeneracy formula of dyonic quarter BPS black holes of
${\mathcal N}=4$ is controlled by the denominator formula of a specific
Borcherds completion of a hyperbolic subalgebra of $\cF$, whose Weyl
group is index $6$ in $W(\cF)$. In more recent work \cite{ChVe} it was
argued that the Weyl group has a more immediate interpretation in this setup
as realizing the so-called attractor flow of a given dyonic solution
to an `immortal dyon' by crossing walls of marginal stability in
moduli space \cite{Be97,Se07,Da07}. The correspondence to the Weyl
group is such that a point $\capX$ inside the forward
light-cone corresponds to dyonic solutions with the given moduli and
charges. As the light-cone is tessellated by the action of the Weyl
group one can move such a point to a standard fundamental domain by a
finite number of Weyl reflections and the endpoint of this motion is
the immortal dyon. The walls between the different cells crossed
during the motion in moduli space are the walls of marginal
stability. It would be interesting to study the degeneracy and the
attractor flow of less supersymmetric solutions and their relation to
the higher rank hyperbolic algebras and modular groups discussed in
this work. The possible relevance of modular forms is also mentioned
in recent work on $E_{10}$ unification and quantum
gravity~\cite{Ganor,DN}. Weyl groups of the so-called hidden symmetry groups
also appear in the context of U-duality
groups~\cite{HuTo94,El97,ObPi98,En07}.\footnote{We stress, however, that
  the Weyl groups discussed here must not be confused with the arithmetic
  duality groups conjectured to be symmetries of
  (compactified) M~theory. For instance, $W(E_{10})$ can be realized 
  as an arithmetic subgroup of $O(1,9)$, whereas the hypothetical arithmetic 
  duality group $E_{10}(\ints)$ is infinitely larger. Independently
  of its possible physical significance, a proper definition of this 
  object at the very least would presuppose properly understanding
  the continuous group $E_{10}$, a goal still beyond reach.}

The organization of the paper is as follows. In Section~\ref{intsec}
we fix the common notation for the four division algebras and describe the
general structure of the hyperbolic root
systems. Section~\ref{weylsec} contains the central general results
about the realization of the simple reflections, the even hyperbolic
Weyl group and the affine and finite Weyl groups for all division
algebras. The commutative division algebras $\reals$ and $\cx$ and the
associated hyperbolic Weyl groups are treated in Section~\ref{comsec},
the algebras admitting a quaternionic realization in
Section~\ref{quatsec} and the case $E_8^{++}$ which requires the
octonions in Section~\ref{octsec}. In appendix~\ref{twapp} we also
present two cases involving complex numbers which are not of
over-extended type but involve a twisted affine extension. 

\vskip10pt
\noindent
{\bf Acknowledgements:}
We would like to thank M.C.N. Cheng, F. Englert, H.R.P. Ferguson,
R.L. Griess, M.~Mazur, W. Nahm, P. Ramond and D. Zagier for discussions 
and correspondence. A.J.F. and A.K. gratefully acknowledge the hospitality 
of the Albert Einstein Institute on various visits, while H.N. thanks
U.L.B. Brussels for hospitality. A.K. is a Research Associate of the 
Fonds de la Recherche Scientifique--FNRS, Belgium.

\end{section}

\begin{section}{Division algebras and hyperbolic root lattices}

We wish to describe the Weyl groups of the various hyperbolic Kac--Moody
algebras as matrix groups which are to be interpreted as modular groups. 
In order to accomplish this, we will describe in this section a lattice
(integral linear combinations of explicit elements) in the Lorentzian space of
Hermitian $(2\times 2)$ matrices over one of the normed division algebras
$\Kk$. This lattice will be shown to be isometric to the root lattice of 
the hyperbolic algebra. In the following section we will show how the Weyl
group acts as on this lattice by matrix conjugation, leading to interesting
matrix groups in the associative cases.

\begin{subsection}{Subrings and integers in division algebras}
\label{intsec}

To begin we review some facts about division algebras, their subrings and
orders in certain algebraic extensions of the rational division algebras. 

Let $\Kk$ be one of the four normed division algebras (over $\reals$):
$\reals$ (the real numbers), 
$\cx$ (the complex numbers), $\quat$ (the Hamilton quaternions) or
$\oct$ (the Cayley octonions). $\Kk$ will also be a topological metric
space with respect to the norm topology. Each algebra $\Kk$ 
has an involution, sending $a\in\Kk$ to $\bar{a}$, generalizing the
complex conjugation in $\cx$, such that the
norm of $a\in\Kk$, $N(a) = a \bar{a}$, is a non-negative real number,
and satisfies  
the composition law $N(ab) = N(a) N(b)$ 
for all $a,b\in\Kk$. For $a\in\reals$, we have $a = \bar{a}$, and for
$a\in\cx$, $\bar{a}$ is  
the usual complex conjugation in $\cx$. One also defines the real part
$a + \bar{a} = 2 Re(a)\in\reals$  
(also sometimes called trace), 
and the imaginary part of $a$ is defined by $a - \bar{a} = 2Im(a)$. 
The real-valued symmetric bilinear form on $\Kk$
\be\label{bilinear}
(a,b) = N(a+b) - N(a) - N(b) = a\bar{b} + b\bar{a}
\ee
is positive definite, giving $\Kk$ the structure of a real Euclidean space.

For each choice of $\Kk$ there is a standard basis such that the structure 
constants are all in $\{-1,0,1\}$, so if one takes only rational linear
combinations of these basis vectors, one gets four rational normed division 
algebras $\Kk_\rats$ (over $\rats$). If $\Ff$ is any subfield of $\reals$ one 
has the normed division algebra $\Kk_\Ff$ (over $\Ff$), consisting of all
$\Ff$-linear combinations of the standard basis elements of $\Kk$, such that  
$\Kk = \Kk_\Ff \otimes_\Ff \reals$.
We will actually only be using a few specific choices for $\Ff$, namely,
$\rats$,  $\rats(\sqrt{2})$, $\rats(\sqrt{3})$ and $\rats(\sqrt{5})$. 

The consideration of $\Kk_\Ff$ is necessary to describe the various finite
root systems which can have angles with cosines involving the listed
square roots. The matrices acting on the root lattices will have entries
from certain {\em orders} $\cO$ within $\Kk_\Ff$. We briefly review
this concept from algebraic number theory (see e.g. \cite{Re75,CoSm03}). The
subfield $\Ff$ of $\reals$, assumed to be a finite extension of $\rats$,
determines a ring $\Ss$ of algebraic integers in $\Ff$, which for our cases
are simply $\ints$, $\ints(\sqrt{2})$, $\ints(\sqrt{3})$ and
$\ints((1+\sqrt{5})/2)$.\footnote{These examples follow the general
  pattern of algebraic integers in $\rats(\sqrt{D})$ which are those
  elements that satisfy a monic minimal polynomial with coefficients
  in $\ints$, that is, $\ints(\sqrt{D})$ or $\ints((1+\sqrt{D})/2)$,
  depending on $D \text{(mod 4)}$.} 
 Note that, as a subring of a field, $\Ss$ is always
commutative. An order  $\cO$ in $\Kk_\Ff$ over $\Ss$ is a subring of 
$\Kk_\Ff$ (non-associative for $\Kk=\oct$) containing $1$ and which is
{\em finitely} generated as an $\Ss$-module (under addition). Furthermore it
has to contain a basis of $\Kk_\Ff$ over $\Ff$, such that 
$\cO\otimes_\Ff \Ff=\Kk_\Ff$. An order is called 
a {\em maximal order} if it is an order not properly contained in another 
order. Orders in rings are a generalization of the integers in the rationals,
and occupy an important role in number theory.

We stress that it is important for an order $\cO$ to be finitely generated and
that this implies that elements $x\in\cO$ satisfy a monic polynomial equation 
\be\label{mp}
x^n + c_{n-1} x^{n-1} + \cdots + c_0 = 0 \qquad
\text{with $c_{n-1}, \dots, c_0 \in\ints$}\,,
\ee
which also can be used for defining orders for the non-associative octonions
\cite{CoSm03}.  
The polynomial (\ref{mp}) arises as follows.
Let $B = \{x_1,\cdots,x_m\}$ be a finite set of generators for $\cO$ over
$\Ss$, so that  
\be
\cO = \left\{\;\sum_{i=1}^m s_i x_i \;\; \vline\;\; s_i\in\Ss\right\}\,.
\ee
Left multiplication by $x\in\cO$, $L_x$, can be represented by a matrix 
$\lp l_{ij}\rp$ such that 
\be
x x_j= \sum_{i=1}^m l_{ij} x_i\quad
\text{and all $l_{ij}\in\Ss$}
\ee
since $\cO$ is a subring and an $\Ss$-module. The characteristic polynomial 
$\det(\lambda I - L_x)$ of the linear transformation $L_x$ is a monic
polynomial of degree $m$ with coefficients in the commutative ring $\Ss$ and
which is satisfied by $L_x$. Since $\Kk$ is a division algebra, the polynomial
is also satisfied by $x$, and so $x$ satisfies a monic polynomial equation with
coefficients in $\Ss$. Multiplying this polynomial by  
the conjugate polynomial, where $\sqrt{D}$ has been replaced by $-\sqrt{D}$
in the coefficients (where $D=2,3,5$), will give a polynomial of degree $2m$
for $x$ with coefficients in $\ints$ as claimed.\footnote{This can also be
 stated as 
 follows: Since $\Ss$ is the ring of algebraic integers over $\Ff$ the norm
 from $\Ff$ to $\rats$ sends $\Ss$ to $\ints$. Applying the norm to the
 equation satisfied by $x$ gives another equation satisfied by $x$, but with
 coefficients in $\ints$.} 
It can happen that the minimal polynomial satisfied by $x$ is of lower degree
but for our purposes it is sufficient to know that some such polynomial
exists. Since in all cases we consider in this paper there is finite group of
units $\Units_\cO$ we know that the ring spanned by them under addition and
multiplication will form an order in the corresponding $\Kk_\Ff$, which we
interpret as integers. 

Although one usually thinks of integers as being `discrete' in some topology, 
note that orders need not be discrete in the norm topology. For example, 
the order $\cO = \ints(\sqrt{2}) = \Ss$ is a dense subring of 
$\Ff = \rats(\sqrt{2})$ where $\Kk = \reals$ (and thus also dense in
$\reals$). Although we will use in some cases such dense
subrings (cf. in particular Section~\ref{quatunits}), the Weyl groups we define
will be discrete groups (in the sense that there exists a fixed $\epsilon >0$
such that the $\epsilon$-neighborhoods of different matrices do not intersect).
As we will explain in Section~\ref{lattsym}, this discreteness can be traced
back to the discreteness of the hyperbolic root lattices. The discreteness
is also important with regard to the properly discontinuous action of
the matrix group on a suitable `upper half-plane' and the existence
of fundamental domains with non-empty interior.\footnote{This can also be
 traced back to the fact that the matrices making up the Weyl group, though
 formally taking entries in some dense order, always are constructed in such
 a way that an underlying discrete order governs the structure, see
 table~\ref{integers}, consistent with the group multiplication laws of the
 matrices.}

For the division algebra $\Kk=\reals$ we take $\Ff = \rats$, $\Ss = \ints$ and
$\cO = \ints$. For $\Kk = \cx$ we use different choices of $\Ff$ and order
$\cO$. Either $\Ff = \rats$, $\Ss = \ints$, 
such that $\Kk_\Ff = \rats(i)$ is the `standard' rational form of $\cx$,
and the  order is the ring of Gaussian integers $\cO = \cG \equiv \ints(i)$;
or we choose $\Ff = \rats(\sqrt{3})$ so that $\Kk_\Ff = \Ff(i)$ is 
an $\Ff$-form of $\cx$ containing the order of Eisenstein integers 
$\cO = \cE \equiv \ints(\root{3}\of{1})$, where  
$\root{3}\of{1} = e^{2\pi i/3} = \frac{-1+i\sqrt{3}}{2}$. 
Each of these is a maximal  order in its algebra
$\Kk_\Ff$. Matters get a little more complicated for $\Kk=\quat$, as is  
to be expected in order to  get the different rank $6$ hyperbolic Weyl
groups.  In one case we choose $\Ff = \rats$, $\Ss = \ints$ and the
(non-maximal) order is the ring of Lipschitz integers $\cO = \Lq$ which are
generated over $\ints$ by the standard basis elements $\{1,i,j,k\}$. The same
rational form contains the maximal order of Hurwitz integers $\cO = \Hq$ which
is generated over $\ints$ by $\{i,j,k,\frac{1+i+j+k}{2}\}$, or by
the elements given below in (\ref{D4roots}).
In another case we add to $\Hq$ the octahedral units (see
Section~\ref{quatunits}), producing 
an order $\cR$ which is generated over $\ints(\sqrt{2})$ by 
$\{a,b,ab,ba\}$ where $a = \frac{1+i}{\sqrt{2}}$ and $b =
\frac{1+j}{\sqrt{2}}$.  
In yet another case we will need to make use of icosian units, and 
appropriate choices of field $\Ff=\rats(\sqrt{5})$, algebraic integers $\Ss$
and order $\cO$, which are described in detail in Section~\ref{quatsec}. 
For the octonions $\Kk=\oct$, and to  
get the Weyl group of $E_{10}$ we will need the octavians which are 
both discrete and a maximal order in the (non-associative) ring of rational
octonions. 

Finally, recall that a unit $\ve\in \cO$ is an element having a multiplicative 
inverse $\nu\in \cO$ in the order, but since $1 = N(\ve\nu) = N(\ve)N(\nu)$ 
and for all the orders we use, $N(\ve),N(\nu)\in\ints$, we find that 
$N(\ve)=\bar{\ve}\ve=\ve\bar{\ve} = 1$. We will use the symbol $\Units_\cO$
to denote the finite group (unless $\Kk = \oct$, where the product is
not associative) of units in $\cO$.

\end{subsection}

\begin{subsection}{Root lattices}

The Kac--Moody algebras we are concerned with are hyperbolic
extensions of finite-dimensional algebras  
associated with discrete subsets of the normed division algebras
$\Kk=\reals,\cx,\quat,\oct$ and, therefore, 
are of ranks $3$, $4$, $6$, and $10$, respectively. In the main text
we consider algebras $\mfg^{++}$ of  
over-extended type \cite{GoOl85,DadBHeSc02} which arise from simple
finite-dimensional Lie algebras  
$\mfg$ by constructing the non-twisted affine extension $\mfg^+$ (also
often denoted by $\mfg^{(1)}$), and connecting the hyperbolic node of
the  
Dynkin diagram by a single line to the affine node. Our methods also
apply to hyperbolic algebras which  
are extensions of twisted affine algebras, and we discuss examples in
Appendix A. But to make the required 
formulas more uniform we restrict ourselves to the non-twisted cases
in the main text.  

We begin with (the Jordan algebra of) all Hermitian $(2\times 2)$ matrices 
over $\Kk$, $\capX = \capX^\dag$,
\be\label{22jordan}
H_2(\Kk) = \left\{\capX = 
 \lp\begin{array}{cc}x^+ & z\\\bar{z} & x^-\end{array}\rp 
 \;\; \vline\;\;
  x^+, \,x^-\in\reals,\,\, z\in \Kk \right\}
\ee 
where $\capX^\dag = \bX^T$ is the conjugate transpose of $\capX$.\footnote{We
  note that all Jordan algebras employed in this paper can also be realized via
  Lorentzian signature Clifford algebras.} 
$H_2(\Kk)$ is equipped with a quadratic
form
\be\label{nrm}
||\capX||^2 = -2\det(\capX)= - 2\big(x^+x^- - z\bar{z}\big)\, 
\ee
and a corresponding symmetric bilinear form
$(\capX,\capY) := \frac12 ( ||\capX + \capY ||^2 - ||\capX ||^2 - 
||\capY ||^2)$ which is Lorentzian. The subspace of matrices 
$\capX\in H_2(\Kk)$ with $x^+ = x^- = 0$ is isomorphic to $\Kk$, and 
the restriction of the bilinear form to this subspace agrees with 
the positive definite form on $\Kk$, making the isomorphism an isometry. 
It is in that Euclidean subspace that we will find root systems
of finite type, and extend them to hyperbolic root systems in $H_2(\Kk)$.

Suppose that we have been able to find a set of simple roots $a_i\in\Kk$, 
$i=1,\ldots,\ell$,  where $\ell = \dim_\reals(\Kk)$, 
with Cartan matrix 
\be\label{Cij}
C= \lp C_{ij} \rp = \lp \frac{2(a_i,a_j)}{(a_j,a_j)}\rp\, .
\ee 
of finite type (the bilinear product is defined in (\ref{bilinear})).
Let $Q = \sum_i \ints a_i \subset \Kk$ be the finite type root 
lattice additively generated by those (finite) simple roots. If necessary, 
we will denote the Lie algebra type by a subscript on $Q$. 
We also find the 
highest root $\theta$ of the finite type root system, and normalize 
the lengths of the simple roots $a_i$ so that $\theta\bth = 1$ always. 
In some twisted affine cases described in Appendix A, we may instead choose
$\theta$ to be the highest short root, but we still assume that 
$\theta\bth = 1$.

In the Jordan algebra $H_2(\Kk)$, depending on the choice of $Q$, we
define the lattice  
\be\label{rtl}
\Lambda = \Lambda(Q) := \left\{\capX = 
 \lp\begin{array}{cc}x^+ & z\\\bar{z} & x^-\end{array}\rp\in H_2(\Kk)
 \;\; \vline\;\;
  x^+, \,x^-\in\ints,\,\, z\in Q  \right\}\,.
\ee
The restriction of the bilinear form to $\Lambda$ makes it a
Lorentzian lattice, which we will identify as the root lattice of a
hyperbolic Kac--Moody Lie algebra. 

Our first step is to identify the hyperbolic simple roots in the lattice 
$\Lambda$ as follows:
\be\label{simroots}
\a_{-1} = \lp\begin{array}{cc}1&0\\0&-1\end{array}\rp,\quad
\a_0 = 
\lp\begin{array}{cc}-1&-\theta\\-\bth&0\end{array}\rp,\quad
\a_i = \lp\begin{array}{cc}0&a_i\\\bar{a}_i&0\end{array}\rp,\ 1\leq
i\leq\ell\,. 
\ee
It follows immediately that
\be
(\a_{-1},\a_{-1}) = 2\;\; ,\quad 
(\a_{-1},\a_0) = -1\;\; ,\quad
(\a_{-1},\a_i) = 0\;\;\hbox{ for }1\leq i\leq\ell\;,
\ee
as well as 
\be
(\a_0,\a_0) &=& 2\theta\bth  = 2 \quad\hbox{and} \nn\\
(\a_i,\a_j) &=& a_i \bar{a}_j + a_j \bar{a}_i = 
(a_i,a_j) \quad \mbox{for $1\leq i,j \leq\ell$.}
\ee
We will give specific choices for the $a_i$ so that the matrix $C$
in (\ref{Cij}) is the Cartan matrix of a finite type root system, 
consistent with the embedding of the corresponding Euclidean 
root lattice $Q$ into $\Lambda$. We will see that choosing $\theta$ to
be the highest root of the finite root system (or the highest short
root) makes $\{\a_0,\a_1,\cdots,\a_\ell\}$ the simple roots of an
untwisted (or a twisted) affine root system. Furthermore we will see that  
including $\a_{-1}$ to get all the simple roots in (\ref{simroots})
gives the simple roots of a hyperbolic root system with Cartan matrix
\be\label{hyperCartan}
[A_{IJ}] := \lp\frac{2(\a_I,\a_J)}{(\a_J,\a_J)}\rp \qquad
\mbox{for $I,J = -1,0,1,\cdots , \ell$.}
\ee
The unit norm condition on $\theta$ is necessary to obtain the single
line between the hyperbolic and affine node. For simply laced algebras
all the $a_i$ can also be chosen as units whereas this is no longer
true for cases where the Dynkin diagram has arrows. 

The hyperbolic Kac--Moody algebras (of over-extended type) of rank 3,4,6 
and 10, respectively, which we obtain by this construction are as follows:
$\cF \equiv A_1^{++}$ (for $\Kk=\reals$), 
$A_2^{++}, C_2^{++}$ and $G_2^{++}$ (for $\Kk=\cx$), 
$A_4^{++}$, $B_4^{++}$, $C_4^{++}$, $D_4^{++}$ and $F_4^{++}$ (for
$\Kk=\quat$),  
and $E_{10} \equiv E_8^{++}$, $D_8^{++}$ and $B_8^{++}$ (for $\Kk
=\oct$).\footnote{The over-extensions of the other finite simple rank $8$ root
 systems $C_8^{++}$ and $A_8^{++}$ do not give rise to hyperbolic reflection
 groups.} 
In all cases, real roots are characterized by
$\det \capX <0$ and imaginary roots by $\det\capX\geq 0$, with null 
roots obeying $\det \capX =0$.

\end{subsection}

\end{section}

\begin{section}{Weyl groups}
\label{weylsec}

In this section we study abstractly the Weyl groups of the 
hyperbolic Kac--Moody algebras identified in the preceding section. 
As we define the root lattice using the Jordan algebra $H_2(\Kk)$ we aim 
to find a description of the Weyl group acting
on this space. Except where explicitly stated otherwise, the results of 
this section apply to all four division algebras, including $\Kk=\oct$. 
Specific features of the four individual division algebras will be 
studied separately in the following sections.

\begin{subsection}{The simple reflections}

The Weyl group associated with the hyperbolic Cartan matrix 
(\ref{hyperCartan}) is the Coxeter group generated by the simple 
reflections 
\be\label{wIdef}
w_I(\capX) = \capX - \frac{2(\capX,\a_I)}{(\a_I,\a_I)} \ \a_I \; ,\quad 
I = -1,0,1,\cdots,\ell\,.
\ee 
These generators are known to satisfy the Coxeter relations, which 
give a complete presentation for the Weyl group. Recall \cite{Hu97} that 
a Coxeter group is defined by the presentation
\be\label{GeneralCoxeter}
\big\langle R_I\ |\ R_I^2 = \id, (R_I R_J)^{m_{IJ}} = \id \,\text{for $I\neq
  J$}\big\rangle \,.
\ee
For Weyl groups only special values may occur for the entries of the matrix 
$M = [m_{IJ}]$, namely $m_{IJ} \in \{2,3,4,6,\infty\}$. More specifically, 
the relations (\ref{GeneralCoxeter}) say that $|w_I| = 2$, that is,
each generator has order $2$, and they give the orders of the products 
of pairs of distinct generators, determined by the entries of the 
Cartan matrix as follows: 
\be\label{CoxeterOrders}
|w_I w_J| = m_{IJ}, \qquad I\neq J \in\{-1,0,1,\cdots,\ell\} ,
\ee
where
\be\label{CoxeterExponents}
m_{IJ} &=&  \left\{\begin{array}{cl}
2 & \text{if $A_{IJ}A_{JI} = 0$}\\
3 & \text{if $A_{IJ}A_{JI} = 1$}\\
4 & \text{if $A_{IJ}A_{JI} = 2$}\\
6 & \text{if $A_{IJ}A_{JI} = 3$}\\
\infty & \text{if $A_{IJ}A_{JI} \geq 4$.}
\end{array}\right.
\ee

For other values of $m_{IJ}$, the Coxeter group defined by 
(\ref{GeneralCoxeter}) cannot be the Weyl group of a Kac--Moody algebra. 
Nevertheless, Coxeter groups with other values of $m_{IJ}$ besides $2,3,4,6,
\infty$ may be of interest in the present context, as such groups may occur
as {\em subgroups} of Weyl groups. A prominent example of such a
non-crystallographic Coxeter group is the group $H_4$, which is a 
subgroup of the Weyl group of $E_8$ \cite{ElSl87,MoPa93,KoKoAl01}.

Of central importance for our analysis is that the action of the
simple hyperbolic Weyl reflections $w_I$ can be rewritten as a matrix
action on $\capX$. 

\begin{thm}\label{Weylprop} Denote by $\bX = (\capX^\dag)^T$ the matrix 
$\capX$ with each entry conjugated but the matrix not transposed, let 
\be\label{ve}
\ve_i=a_i/\sqrt{N(a_i)} \quad \mbox{for $i=1,\cdots , \ell$}
\ee
be the unit norm versions of the simple roots $a_i\in Q$, and let 
$\theta\in Q$ be the highest root of the finite root system with simple 
roots $a_i$, always normalized so that $\theta\bth = 1$. (For simply
laced finite root systems we have $\ve_i = a_i$.) 
Define the $(2\times 2)$ matrices $M_I$, $I = -1,0,1,\cdots,\ell$, by
\be\label{Weyl2}
M_{-1} = \lp\begin{array}{cc}0&1\\1&0\end{array}\rp\,,\quad
M_{0} = \lp\begin{array}{cc}-\theta&1\\0&\E\end{array}\rp\,,\quad
M_{i} = \lp\begin{array}{cc}\ve_i&0\\0&-\bar{\ve_i}\end{array}\rp\,,
\ee
$i = 1,\cdots,\ell$. Then the simple reflections (\ref{wIdef}) can be
written as 
\be\label{MWeyl}
w_I(\capX) = M_I \bar{\capX}M_I^\dagger
\quad,\ I = -1,0,1,\ldots,\ell.
\ee
\end{thm}

\noindent
{\bf{Proof:}} In the non-associative octonionic case one needs to check 
that the expression on the right side of (\ref{MWeyl}) is well defined 
for the matrices in (\ref{Weyl2}) without placing parentheses. Since
each matrix $M_I$  
involves only one non-real octonion, this follows immediately from 
the alternativity property of $\oct$. (See the section on octonions
below for more details.)  

Now we will check directly that the two formulas 
(\ref{wIdef}) and (\ref{MWeyl}) for $w_I$ agree. First, note that
\be\label{MXMI}
(\capX,\a_{-1}) = x^+ - x^-, \quad (\capX,\a_0) = -z \bth - \theta
\bar{z} + x^-, \quad 
(\capX,\a_i) = z \bar{a}_i + a_i \bar{z}
\ee
for $1\leq i\leq \ell$, so
\be
w_{-1}(\capX) &=& \capX - \frac{2(\capX,\a_{-1})}{(\a_{-1},\a_{-1})} \ 
\a_{-1} = \capX - (x^+ - x^-) \a_{-1} = 
  \lp\begin{array}{cc}x^- & z\\\bar{z} & x^+\end{array}\rp \, ,\nn\\
w_0(\capX) &=& \capX - \frac{2(\capX,\a_0)}{(\a_0,\a_0)} \ \a_0 = 
\capX - (-z \bth - \theta \bar{z} + x^-) \a_0    \nn\\
&=& \lp\begin{array}{cc}(x^+ -
z\bth-\theta\bar{z}+x^-)&(z-z\bth\theta-\theta\bar{z}\theta+x^-\theta)
\\ 
(\bar{z}-\bth z\bth-\bth\theta\bar{z}+\bth x^-) & x^-\end{array}\rp \, ,\nn\\
w_i(\capX) &=& \capX - \frac{2(\capX,\a_i)}{(\a_i,\a_i)} \ \a_i = 
\capX - \frac{z \bar{a}_i + a_i \bar{z}}{N(a_i)} \a_i \nn\\   
&=&  \lp\begin{array}{cc}x^+ &  -\ve_i \bar{z} \ve_i\\ 
- \bve_i z\bve_i & x^-\end{array}\rp .
\ee
where we used the definition (\ref{ve}) in the last equation.
Compare these with 
\be\label{MXMII}
M_{-1} \bX M_{-1}^\dag &=& \lp\begin{array}{cc}0&1\\1&0\end{array}\rp\ 
\lp\begin{array}{cc}x^+ & \bar{z}\\ z & x^-\end{array}\rp\ 
\lp\begin{array}{cc}0&1\\1&0\end{array}\rp
= \lp\begin{array}{cc}x^- & z\\\bar{z} & x^+\end{array}\rp  \, , \nn\\
M_0 \bX M_0^\dag &=& \lp\begin{array}{cc}-\theta &1\\0&\bth\end{array}\rp\ 
\lp\begin{array}{cc}x^+ & \bar{z}\\ z & x^-\end{array}\rp\ 
\lp\begin{array}{cc}-\bth &0\\1&\theta\end{array}\rp    \nn\\
&=& \lp\begin{array}{cc}(\theta\bth x^+ - z\bth - \theta\bar{z}+x^-) &
(-\theta\bar{z}\theta+x^-\theta) \\ 
(-\bth z\bth + \bth x^-) & \bth\theta x^-\end{array}\rp \, , \nn\\
M_i \bX M_i^\dag &=& 
\lp\begin{array}{cc}\ve_i&0\\0&-\bar{\ve_i}\end{array}\rp\ 
\lp\begin{array}{cc}x^+ & \bar{z}\\ z & x^-\end{array}\rp\ 
\lp\begin{array}{cc}\bar{\ve_i}&0\\0&-\ve_i\end{array}\rp \nn\\
&=& \lp\begin{array}{cc}\ve_i\bar{\ve}_i x^+ & -\ve_i \bar{z} \ve_i \\
- \bar{\ve_i} z \bar{\ve_i}& \bar{\ve_i}\ve_i x^-\end{array}\rp.
\ee
Hence $w_{-1}(\capX) = M_{-1} \bX M_{-1}^\dag$, 
$w_0(\capX) = M_0 \bX M_0^\dag$, since $\theta\bth = 1$, and 
$w_i(\capX) = M_i \bX M_i^\dag$ since $\ve_i\bve_i =1$.  The explicit
expressions (\ref{MXMI}) and (\ref{MXMII}) are manifestly well defined
for octonions without placing parentheses. $\square$ 
\medskip

Our aim in this paper is to study the group generated by the simple 
reflections (\ref{MWeyl}) as an extension of a matrix group. This 
extension may always be realized in the associative cases ($\Kk \neq \oct$)
as an extension of a matrix group by a small finite group. Since the 
full Weyl group $W$ is a semi-direct product of the even part $W^+$ with 
$\Ztwo = \langle w_{-1}\rangle$, we will be satisfied to understand just 
the even Weyl group $W^+$ as an extension of a matrix group.

\end{subsection}

\begin{subsection}{Even part of the Weyl group}

The formula (\ref{MWeyl}) for the simple reflections involves complex
conjugation of $\capX$. Therefore all {\em even} elements $s\in W^+\subset W$
can  be represented {\em without complex conjugation of $\capX$}, and it 
turns out to be simpler to study the even Weyl group $W^+$ in many cases.

The even Weyl group $W^+$ is an index 2 normal subgroup of $W$ and consists of
those elements which can be expressed 
as the product of an even number of simple reflections. It is generated by the 
following list of $\ell + 1$ double reflections:
\be\label{evenweylgen}
s_{0} = w_{-1}w_0\,,\quad s_i =w_{-1}w_i\quad(i=1,\ldots,\ell)\,.
\ee
Of course, this is not a unique set of generating elements. From the Coxeter
relations  (\ref{GeneralCoxeter})--(\ref{CoxeterExponents}) satisfied by the
simple reflections, $w_I$, these even elements satisfy the relations
\be
s_0^3 = \id\,,\quad s_i^2=\id\,\,\,\text{  for $i\neq 0$}
\ee
and
\be\label{evenweylrel}
(s_i^{-1}s_j)^{m_{ij}} = \id \,\,\,\text{for $i\neq j$\, and\,
  $i,j=0,1,\ldots,\ell$}, 
\ee
where $m_{ij}$ is given as before in (\ref{CoxeterExponents}). 

In analogy with Theorem~\ref{Weylprop} we have

\begin{thm}\label{EvenW} Define the matrices
\be\label{W+}
S_0 =
   \lp\begin{array}{cc}0&\theta\\-\E&1\end{array}\rp\,,\quad
S_i =
   \lp\begin{array}{cc}0&-\ve_i\\\bar{\ve}_i&0 \end{array}\rp\,.
\ee
Then the generating double reflections (\ref{evenweylgen}) acting on
$\capX\in\Lambda$ can be written for all $\Kk$ as 
\be\label{EvenWMat}
s_I(\capX) = S_I \capX S_I^\dagger \quad, \, I=0,1,\ldots,\ell\,.
\ee
\end{thm}

\noindent{\bf Proof:} Follows by direct computation as in the proof of
Theorem~\ref{Weylprop}. $\square$ 
\medskip

An important corollary is
\begin{cor}\label{RCHeven}
For the associative division algebras $\Kk=\reals,\cx,\quat$,
{\em all} elements $s\in W^+$ can be realized by a matrix action
according to
\be\label{evenmat}
s(\capX) = S\capX S^\dagger\,,
\ee
where $S=S_{i_1}\cdots S_{i_n}$ if $s=s_{i_1}\cdots s_{i_n}\in W^+$ in terms
of the generating elements (\ref{evenweylgen}). 
\end{cor}

\noindent{\bf Proof:} 
The iterated action of two even Weyl transformations is given by the
associative product  of matrices (for $\Kk=\reals,\cx,\quat$) 
\be
(s_1 s_2)(\capX) = S_1(S_2\capX S_2^\dagger)S_1^\dagger =
(S_1S_2)\capX(S_1S_2)^\dagger\, 
\ee
and has an obvious extension to arbitrary words in the even Weyl group by
associativity. $\square$ 
\medskip

For $\Kk=\oct$, the formula (\ref{evenmat}) no longer holds (even though we
have for any octonionic matrices $(S_1S_2)^\dagger = S_2^\dagger
S_1^\dagger$). This can be seen most easily in the continuous case by a
dimension count as will be discussed in Section~\ref{octsec} where we
collect our results specific to the octonionic case.

Let us remark that in the commutative case (\ref{EvenWMat}) follows
immediately by acting with 
two successive simple Weyl reflections, say $w_J$ and $w_I$; the
effect on $\capX$ is 
\be
w_I(w_J(\capX)) = M_I
\overline{(M_J\bar{\capX}M_J^\dagger)}M_I^\dagger
= S_{IJ} \capX S_{IJ}^\dagger  \; ,
\ee
where $S_{IJ} \equiv M_I \bM_J$. In this notation the matrices (\ref{W+}) are
$S_i \equiv S_{-1\,i}$ for $i=0,1,\ldots,\ell$. However, in the
non-commutative cases one has to be more careful because quaternionic and
octonionic conjugation also reverses the order of factors inside a product,
such that, in general  
\be\label{noncprob}
\overline{M_J\bar{\capX}M_J^\dagger}\ne \bar{M}_J \capX M_J^T
\ee
and $M_I$ and $\bM_J$ do not obviously combine into a matrix $S_{IJ}$ which
acts by conjugation as in (\ref{evenmat}).  Therefore it is crucial that
Theorem~\ref{EvenW} applies to {\em all} division algebras. 

In the associative cases we obtain from Corollary~\ref{RCHeven} that $W^+$ is
isomorphic to a matrix group generated by (\ref{W+}). All elements $s\in W^+$
act by  invertible matrices and we therefore obtain subgroups of $GL_2(\Kk)$.
Furthermore, they all act by matrix conjugation and therefore a matrix and its
negative have the same action. Other scalar matrices $\ve \id \in GL_2(\Kk)$,
for $\ve$ a central unit in $\Kk$, would act trivially, so we should be
finding $W^+$ isomorphic to a subgroup of $PGL_2(\Kk)$. It is important to
note that for general $(2\times 2)$ matrices over $\Kk$ the determinant is not
well defined unless $\Kk$ is commutative. In the cases $\Kk=\reals$ and $\Kk =
\cx$, the matrices $M_I$ satisfy $\det M_I = -1$, so products of two such
fundamental Weyl reflections have determinant $+1$ and hence $W^+$ is a
subgroup of $PSL_2(\Kk)$ for commutative $\Kk$. With a suitable definition of
$PSL_2(\Kk)$ this statement is also true for non-commutative $\Kk$. The
necessary refinements required for  $PSL_2(\quat)$ and $PSL_2(\oct)$ will be
presented in the relevant sections. For all $\Kk$ the
$P$ in $PSL$ means that only the quotient by $\{\id, -\id\}$ has been
taken. We will argue that the even Weyl groups
constitute interesting discrete `modular' subgroups of $PSL_2(\Kk)$.

\end{subsection}

\begin{subsection}{Finite and affine Weyl subgroups}
\label{faff}

By specialization, the construction given above also yields matrix 
representations of the finite and the affine Weyl subgroups contained in 
the respective hyperbolic algebras. We first note that the action of 
the {\em finite} Weyl group $W_{\rm fin} \equiv W(\mfg)$ on the root 
lattice of the finite subalgebra $\mfg$ is obtained as a special case 
of (\ref{MWeyl}) by setting $x^\pm = 0$ and restricting indices to
$I\equiv i=1,\cdots,\ell$; as follows immediately from the last formula 
in (\ref{MXMII}), the simple reflections are thus realized on any lattice 
vector $z$ via
\be\label{MWeylf}
w_i(z) = - \ve_i \bar{z}\ve_i \qquad\qquad \mbox{(for $z\in Q\subset \Kk$)}. 
\ee
The {\em same} transformation on $z$ is obtained by matrix conjugation 
(\ref{evenmat}) with the even element $S_i$ with $i =1,\ldots,\ell$, i.e. 
{\em without} complex conjugation of $\capX$. Similarly, it follows from 
(\ref{evenweylgen}) that 
\be
s_is_j = w_iw_j\quad \text{for}\,\,i,j=1,\ldots,\ell\,.
\ee
so that, in terms of the matrix representation (\ref{W+}) for $\Kk\neq \oct$ 
we obtain 
\be\label{SiSj}
S_iS_j =
\lp\begin{array}{cc}-\ve_i\bar{\ve}_j&0\\
 0&-\bar{\ve}_i\ve_j\end{array}\rp
\ee
whence the even part $W^+(\mfg)$ of the finite Weyl group acts by diagonal
matrices. Defining $u_{ij}=\ve_i\bar{\ve}_j$ and $v_{ij}=\bar{\ve}_i\ve_j$ we
deduce $s_i (s_j (z)) =  
u_{ij} z \bar{v}_{ij} = u_{ij} z v_{ji}$ for $z\in Q$. To summarize: the even
and odd parts of  
the finite Weyl group, respectively, act by purely diagonal or purely 
off-diagonal matrices for $\Kk=\reals,\cx,\quat$.~\footnote{\label{csnot} For
 $\Kk=\quat$  
the same action can be alternatively written in terms of pairs $[l,r]$ 
and $*[l,r]$ of unit quaternions, cf. \cite{CoSm03} p.~42. The relation to our
notation is as follows:
\be
[l,r] \longleftrightarrow 
  S=\lp\begin{array}{cc}\bar{l}&0\\
           0&\bar{r}\end{array}\rp\,,\quad\quad
*[l,r] \longleftrightarrow 
  S=\lp\begin{array}{cc}0&\bar{l}\\
           \bar{r}&0\end{array}\rp\,.\nn
\ee}
Because of the extra $w_{-1}$ contained in the definition of $S_i$ the 
action of the odd part is only correct on the subspace with $x^\pm=0$; on the
full $H_2(\Kk)$ the odd parts act with an additional interchange of $x^+$ and
$x^-$. 

The affine subalgebra $\mfg^+$ is characterized by all roots with 
$\det\capX \leq 0$ which are of the form
\be
\capX = \lp\begin{array}{cc}m& z\\  \bar{z}& 0\end{array}\rp \qquad
\mbox{for $z\in Q$ and $m\in\ints$}\,.
\ee
Its associated affine Weyl group $W_{\rm aff}$ is well known to be isomorphic 
to a semi-direct product of the finite Weyl group and an abelian group of
(affine) translations $\cT$, such that $W_{\rm aff} = \cT\rtimes W_{\rm fin}$. 
The latter is generated by elements of the form $w t_\theta w^{-1}=
t_{w(\theta)}$ 
where $w\in W_{\rm fin}$ and the relevant affine translation is
\be\label{afftrl}
t_\theta = w_0 w_\theta \; = 
   \lp\begin{array}{cc} 1&\theta\\ 0& 1\end{array}\rp \,.
\ee
Here, $w_\theta$ is the reflection about the highest root $\theta$ and
this is the correct expression for the  affine translation for all the
algebras we consider in this paper. It 
is straightforward to check that this matrix indeed generates 
translations since
\be
t_\theta \lp\begin{array}{cc}m& z\\  \bar{z}& 0\end{array}\rp t_\theta^\dagger
= \lp\begin{array}{cc}m + \theta \bar{z} + z\bth& z\\ \bar{z}& 0\end{array}\rp 
\ee
This statement also holds for the octonionic case.
We note that the interesting $S$-type transformations in the Weyl
groups are then solely due to the hyperbolic extension.

If $\Kk$ is associative, and once the finite Weyl group has been
identified in terms of diagonal and off-diagonal matrices, the full
even hyperbolic Weyl group is obtained by adjoining the affine
Weyl transformation (\ref{afftrl}) to the set of both diagonal and off-diagonal
matrices as generating set. 

\end{subsection}

\begin{subsection}{Lattice symmetries and Weyl groups}
\label{lattsym}

In the associative cases, there are general necessary constraints on the 
structure of the matrices
\be
S = \lp \begin{array}{cc} a&b\\c&d\end{array}\rp\,,
\ee
which represent transformations of the even Weyl group acting on the
hyperbolic root lattice $\Lambda$ of (\ref{rtl}). These follow from
the fact that 
the transformed $\capX'=S\capX S^\dagger$ again has to lie in the root
lattice $\Lambda$ and should have the same norm as $\capX$. 
Working out the matrix product one finds the
following set of conditions for the transformed $\capX'$ to lie in the
root lattice:
\be\label{rtsym}
a\bar{a},\,b\bar{b},\,c\bar{c},\,d\bar{d}&\in& \ints\,,\nonumber\\
c\bar{a},\,d\bar{b}&\in& Q\,,\nonumber\\
a a_i \bar{b} + b \bar{a}_i \bar{a},\, c a_i
\bar{d}+d\bar{a}_i\bar{c}&\in& \ints\;\quad\text{for
 $i=1,\ldots,\ell$}\,,\nonumber\\ 
a a_i\bar{d} + b\bar{a}_i\bar{c}&\in& Q\;\quad\text{for
 $i=1,\ldots,\ell$}\,.
\ee
Here, $Q$ is the finite root lattice with simple basis vectors
$a_i$. It is the collection of the conditions (\ref{rtsym}) that will 
turn the matrices $S$ into a {\em discrete subgroup} of the group of 
$(2\times 2)$ matrices over $\Kk=\reals,\cx,\quat$ (even if the 
relevant rings inside $\Kk$ to which the matrix entries belong are 
not discrete). The norm preservation requirement will lead to additional 
determinant-type constraints. The conditions (\ref{rtsym}) are only 
necessary but not sufficient since they also 
allow for solutions which correspond to lattice symmetries of $Q$
which are not elements of the Weyl group. This happens when the finite
Dynkin diagram admits outer automorphisms, and we will be concerned
with finding manageable conditions which eliminate these.

\end{subsection}

\end{section}

\begin{section}{Commutative cases}
\label{comsec}

We first discuss the commutative cases $\Kk=\reals$ and $\Kk=\cx$
where one has the usual definition of the determinant. 

\begin{subsection}{$\Kk= \reals$, type $A_1$}

For this case we recover the results of \cite{FeFr83} for the rank 3 
hyperbolic algebra $\cF = A_1^{++}$, where $W(\cF) = PGL_2(\ints)$. 
The root system of type $A_1$ is shown in fig.~\ref{a1latt}, and the
Dynkin diagram of $\cF$ is shown in fig.~\ref{a1pp}. 
A simplification here is that we do not have to worry about
conjugation; it is for this reason that an isomorphism with a matrix
group exists for the full Weyl group $W$, rather than only its even 
subgroup $W^+$. We have the simple root $a_1=1$, which is identical to 
the highest root $\theta =1$, so that the simple roots of the
hyperbolic algebra  
are represented by the three matrices
\be\label{simrootsF}
\a_{-1} = \lp\begin{array}{cc}1&0\\0&-1\end{array}\rp\quad,\quad
\a_0 = \lp\begin{array}{cc}-1&-1\\-1&0\end{array}\rp\quad,\quad
\a_1 = \lp\begin{array}{cc}0&1\\ 1&0\end{array}\rp\,.
\ee
The corresponding matrices $M_I$ representing the simple Weyl group
generators are 
\be\label{MFF}
M_{-1} = \lp\begin{array}{cc}0&1\\1&0\end{array}\rp\quad,\quad
M_0 = \lp\begin{array}{cc}-1& 1\\ 0&1\end{array}\rp\quad,\quad
M_1 = \lp\begin{array}{cc}1& 0\\ 0& -1\end{array}\rp\,.
\ee
The even part of the Weyl group is thus generated by the following
two matrices (conjugation can be omitted)
\be
S := M_{-1} M_1 =
\lp\begin{array}{cc}0&-1\\1&0\end{array}\rp\;,\quad
T := M_0 M_1 =
\lp\begin{array}{cc}-1& -1\\ 0 & -1 \end{array}\rp \cong
\lp\begin{array}{cc} 1&  1\\ 0 & 1 \end{array}\rp
\ee
implying that $W^+ (\cF) \cong PSL_2 (\ints)$. The full Weyl group in
this case can be obtained by adjoining the matrix $M_{-1}$ of
determinant $-1$ from which one recovers $W (\cF) \cong PGL_2
(\ints)$ \cite{FeFr83}.

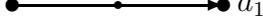
\begin{figure}
\centering
\begin{picture}(50,90)
\put(65,28){$a_1$}
\put(-20,30){\line(1,0){80}}
\put(-20,30){\circle*{5}}
\put(20,30){\circle*{3}}  
\put(60,30){\circle*{5}}
\thicklines
\put(20,30){\vector(1,0){40}}
\end{picture}
\caption{\label{a1latt}\sl The $A_1$ root system with simple root labeled. 
  The root lattice is the lattice of (rational) integers $\ints$.}
\end{figure}

\begin{figure}
\centering
\begin{picture}(100,20)
\thicklines
\multiput(10,10)(30,0){3}{\circle*{8}}
\put(14,10){\line(1,0){22}}
\multiput(40,7)(0,2){4}{\line(1,0){30}}
\put(4,-5){$-1$}
\put(37,-5){$0$}
\put(67,-5){$1$}
\end{picture}
\caption{\label{a1pp}\sl Dynkin diagram of $\cF \equiv A_1^{++}$
 with numbering  of nodes.}
\end{figure}

\end{subsection}

\begin{subsection}{$\Kk=\cx$, simple reflections}

For $\Kk=\cx$ there are different choices of simple finite root systems
which we discuss separately in the following sections, corresponding to 
the root lattices $A_2,B_2\cong C_2, G_2$. Here, we collect some 
common features of all cases. However, the overextension of the 
$B_2\cong C_2$ root system leads to the hyperbolic algebra $C_2^{++}$
if $\theta$ is the highest root. 

In all cases in this subsection, the matrices giving rise to the
simple reflections are  
\begin{align}
M_{-1} &=\lp\begin{array}{cc}0&1\\1&0\end{array}\rp\,,& 
M_{0}&=\lp\begin{array}{cc}-\theta&1\\0&\E\end{array}\rp\,,&\nn\\
M_{1} &=\lp\begin{array}{cc}1&0\\0&-1\end{array}\rp\,,& 
M_{2} &=\lp\begin{array}{cc} \ve_2&0\\0& - \bve_2\end{array}\rp\,.&
\end{align}
where we used the fact that for $\Kk=\cx$ it is always possible
to choose $\ve_1 =1$ (as for $\Kk=\reals$). The even Weyl group 
$W^+$ is thus generated by the elements
\be\label{SC}
S_0 &=& M_{-1} \bM_0 
  = \lp\begin{array}{cc}0&\theta\\-\E&1\end{array}\rp\,,\quad
S_1 = M_{-1} \bM_1
  = \lp\begin{array}{cc}0&-1\\1&0\end{array}\rp\,,\nn\\
S_2 &=& M_{-1} \bM_2 = 
   \lp\begin{array}{cc}0&-\ve_2\\ \bve_2&0\end{array}\rp\,
\ee
and we note that $S_2^2 = -\id$ which acts as the identity on $\capX$ in
agreement with the Coxeter relation $s_2^2=\id$. 
We repeat that we will always normalize $\theta\bth=1$. In the above 
form the main difference between the three algebras is encoded 
in  $\theta$ and $\ve_2$. The three algebras are then distinguished 
simply by the multiplicative order of these numbers, which will be 
specified below for each case.

We note that the even parts of the finite Weyl groups $W^+(A_2)$, 
$W^+(B_2)\cong W^+(C_2)$ and $W^+(G_2)$ are 
cyclic groups of orders $3$, $4$, and $6$, respectively. This follows 
also since their generating elements are rotations in the plane and 
so one obtains finite subgroups of the abelian group $SO(2)$. From 
this point of view the non-abelian nature of the full finite Weyl group
arises because of a single reflection realized as complex conjugation,   
which shows that the full finite Weyl group is a dihedral group.

In order to determine the even Weyl group $W^+$ for the hyperbolic
algebras $A_2^{++}$ and $G_2^{++}$ the following result will be useful. 

\begin{prop}\label{sl2prop}
Let $\cO$ be a discrete Euclidean ring in $\cx$, that is, a discrete
additive group,  
closed under multiplication, satisfying the Euclidean
algorithm. Furthermore, assume that  
$\cO = \{m+n\theta^2\in\cx\ |\ m,n\in\ints\}$ with $\theta\in \cO$
a unit, and that  
all units in $\cO$ have norm $1$. 
Then $SL_2(\cO)$, the group of $(2\times 2)$ matrices 
with entries from $\cO$ and with determinant $1$, is generated by
\be\label{sl2cx}
\left(\begin{array}{cc}1&1\\0&1   \end{array}\right)\,,
\left(\begin{array}{cc}\ve&0\\0&\bar\ve   \end{array}\right)\,,
\left(\begin{array}{cc}0&1\\-1&0   \end{array}\right)\,,
\ee
where $\ve\in\Units_\cO$ runs through all units of $\cO$.
\end{prop}

\noindent{\bf{Proof:}} 
We prove the proposition by using arguments from \cite{Kr85}. Let
$\Delta$ denote the group  
generated by the matrices in (\ref{sl2cx}). These matrices belong to
$SL_2(\cO)$ and we want  
to show that $\Delta=SL_2(\cO)$. First, we claim that $\Delta$
contains all `translation' matrices 
\be\label{trmatscx}
\left(\begin{array}{cc}1&a\\0&1   \end{array}\right)\quad \text{for}\,
a\in  \cO. 
\ee 
To see this, note that 
\be\label{trmat}
\left(\begin{array}{cc}\theta&0\\0&\bth \end{array}\right)
\left(\begin{array}{cc}1&1\\0&1   \end{array}\right)
\left(\begin{array}{cc}\bth&0\\0&\theta \end{array}\right)
= \left(\begin{array}{cc}1&\theta^2\\0&1   \end{array}\right)\,.
\ee
We used that the unit $\theta\in \cO$ satisfies $\theta\bth =
1$. Since the translation matrices form  
an abelian group with addition of the upper right corner entry, and
$1$ and $\theta^2$ are an integral  
basis of $\cO$, for any $m,n\in\ints$, we have
\be
\left(\begin{array}{cc}1&1\\0&1   \end{array}\right)^m 
\left(\begin{array}{cc}1&\theta^2\\0&1   \end{array}\right)^n
= \left(\begin{array}{cc}1&m+n\theta^2\\0&1   \end{array}\right)\,,
\ee
so all matrices (\ref{trmatscx}) are in $\Delta$. We also get that all
matrices  
\be
\left(\begin{array}{cc}0&1\\-1&0 \end{array}\right)
\left(\begin{array}{cc}1&-a\\0&1 \end{array}\right)
\left(\begin{array}{cc}0&-1\\1&0 \end{array}\right)
= \left(\begin{array}{cc}1&0\\a&1 \end{array}\right)
\ee
are in $\Delta$. To complete the proof, choose any matrix 
$A\in SL_2(\cO)$ and consider the set of norms
\be
\cN =\big\{ N(b_{ij}) \ |\ [b_{ij}] = B = UAV \text{ for some
}U,V\in\Delta\big\}\backslash\{0\}. 
\ee
Because the ring $\cO$ is discrete, and the norm is positive
definite, any set of non-zero norms has  
a least element. Then $\cN$ contains a non-zero minimum element
$N(b)$, and suppose 
\be 
B = \left(\begin{array}{cc}b_{11}&b_{12}\\b_{21}&b_{22}\end{array}\right)
\ee
is a matrix for which this occurs. Multiplying $B$ on the left and/or
right by the third generating matrix  
in (\ref{sl2cx}) (the rotation matrix), we can move any of its entries
into the upper left corner, so we may 
assume that the minimum occurs with $N(b)=N(b_{11})$. Now, using the
Euclidean algorithm, we may write  
\be
b_{12} &=& q_1 b_{11} + r_1 \quad \text{with $0\le N(r_1) < N(b_{11})$}\nn\\
b_{21} &=& q_2 b_{11} + r_2 \quad \text{with $0\le N(r_2)<N(b_{11})$}
\ee
for $q_1,q_2,r_1,r_2\in \cO$. Then we have
\be
B \left(\begin{array}{cc}1&-q_1\\0&1\end{array}\right) = 
\left(\begin{array}{cc}b_{11}&r_1\\b_{21}&b_{22}-q_1 b_{21}\end{array}\right)
\ee
is a matrix of the form $UAV$ for $U,V\in\Delta$, so the norms of its
entries cannot be less than the 
minimal value $N(b_{11})$. This forces $r_1 = 0$, so $b_{12} = q_1
b_{11}$. Similarly, we have 
\be
\left(\begin{array}{cc}1&0\\-q_2&1\end{array}\right) B =
\left(\begin{array}{cc}b_{11}&b_{12}\\r_2&b_{22}-q_2 b_{12}\end{array}\right)
\ee
is a matrix of the form $UAV$ for $U,V\in\Delta$, so the norms of its
entries cannot be less than the 
minimal value $N(b_{11})$. This forces $r_2 = 0$, so $b_{21} = q_2
b_{11}$. Finally, we see that 
\be
\left(\begin{array}{cc}1&0\\-q_2&1\end{array}\right) 
B \left(\begin{array}{cc}1&-q_1\\0&1\end{array}\right)
= \left(\begin{array}{cc}b_{11}&0\\0&b'_{22}\end{array}\right) = B' \,,
\ee
where $b'_{22}= b_{22}-q_1q_2b_{11}$, is a diagonal matrix of the form
$UAV$ for $U,V\in\Delta$, with determinant $b_{11}b'_{22} = 1$.  
Therefore, $b_{11}$ must be a unit in $\cO$ and $b'_{22} =
b_{11}^{-1} = \overline{b_{11}}$  
so that $B'\in\Delta$. All operations transforming $A$ into $B'\in\Delta$ 
were performed using matrices from $\Delta$, so we conclude that 
$A\in\Delta$, completing the proof. $\square$\medskip

We remark that this proposition does not apply to the Gaussian integers since
all Gaussian units square to real numbers and one therefore cannot generate
the whole ring from $1$ and $\theta^2$ for any unit $\theta$.

\end{subsection}

\begin{subsection}{$\Kk=\cx$, type $A_2$}

The first choice of integers we consider is the case of type $A_2$, which 
is simply laced. The simple roots can therefore be chosen as units 
$a_i=\ve_i$ and we take them to be
\be\label{typeA2}
\ve_1 = a_1 =1 \;\; , \; \ve_2 = a_2  = \frac{-1 + i\sqrt{3}}{2} 
  \;\; , \;\; \theta = - \bar{\ve}_2=\frac{1+i\sqrt{3}}{2} \nn\,,
\ee
where also the highest root $\theta$ has been given.
The $A_2$ root lattice is spanned by integral linear combinations of the 
simple roots, and they form the order of `Eisenstein integers' $\cE$. 
The $A_2$ root system is depicted in fig.~\ref{a2roots}.

\begin{figure}
\centering
\begin{picture}(90,90)
\put(40,0){\line(-2,3){40}}
\put(65,28){$a_1$}
\put(-12,65){$a_2$}
\put(50,62){$\theta$}
\put(-20,30){\line(1,0){80}}
\put(0,0){\line(2,3){40}}
\put(20,30){\circle*{3}}  
\put(-20,30){\circle*{5}}
\put(0,0){\circle*{5}}
\put(40,0){\circle*{5}}
\put(60,30){\circle*{5}}
\put(0,60){\circle*{5}}
\put(40,60){\circle*{5}}
\thicklines
\put(20,30){\vector(-2,3){20}}
\put(20,30){\vector(1,0){40}}
\end{picture}
\caption{\label{a2roots}\sl The $A_2$ root system with simple roots labeled. 
  The root lattice is the ring of Eisenstein integers.}
\end{figure}
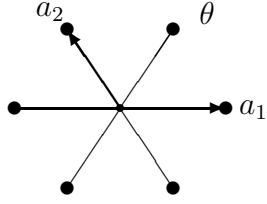

The hyperbolic algebra $A_2^{++}$ has the Dynkin diagram shown 
in fig.~\ref{a2pp}. The following choice of simple roots in $\Lambda$ 
provides us with that diagram:
\begin{align}
\a_{-1}&=\lp\begin{array}{cc}1&0\\0&-1\end{array}\rp\,,&
\a_0&=\lp\begin{array}{cc}-1&-\theta\\-\bar\theta&0\end{array}\rp\,,&\nn\\
\a_{1}&=\lp\begin{array}{cc}0&1\\1&0\end{array}\rp\,,&
\a_{2}&=\lp\begin{array}{cc}0&-\bar\theta\\-\theta&0\end{array}\rp\,,&
\end{align}
where we used $\ve_2=-\bar{\theta}$.
All roots have equal length in this case, so the hyperbolic extension
$A_2^{++}$ is also simply laced. We have

\begin{prop}\label{A2++Weyl} The even part of the Weyl group in this case is
\be\label{W+A2}
W^+(A_2^{++}) \cong PSL_2(\cE)\,
\ee
where $PSL_2(\cE)$ denotes the `Eisenstein modular subgroup' of $PSL_2(\cx)$ 
obtained by restricting all entries to be Eisenstein integers.
\end{prop}

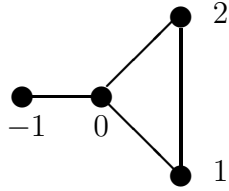
\begin{figure}
\centering
\begin{picture}(110,80)
\thicklines
\multiput(10,40)(30,0){2}{\circle*{8}}
\put(70,70){\circle*{8}}
\put(70,10){\circle*{8}}
\put(14,40){\line(1,0){22}}
\put(42,44){\line(1,1){28}}
\put(42,38){\line(1,-1){28}}
\put(70,66){\line(0,-1){52}}
\put(4,25){$-1$}
\put(37,25){$0$}
\put(82,68){$2$}
\put(82,8){$1$}
\end{picture}
\caption{\label{a2pp}\sl Dynkin diagram of $A_2^{++}$ with numbering
 of nodes.}
\end{figure}

\noindent{\bf{Proof:}}
The statement (\ref{W+A2}) is an immediate corollary of Proposition~\ref{sl2prop},
given that $\cE$ satisfies the Euclidean algorithm, and all
generating matrices in (\ref{sl2cx}) can be obtained within 
the even part of the Weyl group of $A_2^{++}$.  This is true by inspection 
of the matrices (\ref{SC}): the rotation matrix is just $S_1^{-1}$, and
all diagonal matrices are obtained from powers of 
\be
\cb:= S_1S_2 = \lp\begin{array}{cc}\theta &0\\0&\bth\end{array}\rp\,.
\ee
Finally the translation matrix is obtained from
\be
(S_2 S_0) (S_1 S_2)^{-1} =
\lp\begin{array}{cc} -\bth^2 & \bth \\0&-\theta^2\end{array}\rp\cdot
\lp\begin{array}{cc} \bth & 0 \\0&\theta\end{array}\rp =
\lp\begin{array}{cc}1& 1 \\0&1\end{array}\rp\,
\ee
where we used $\theta^3= -1$. 
Thus $S_0$, $S_1$ and $S_2$ generate a group isomorphic to $SL_2(\cE)$. 
Since in the action, the normal subgroup $\{\id, -\id\}$ acts trivially,  
the action of the even Weyl group on $\capX$ is that of the quotient 
$PSL_2(\cE)$. $\square$

\end{subsection}

\begin{subsection}{$\Kk=\cx$, type $C_2$}

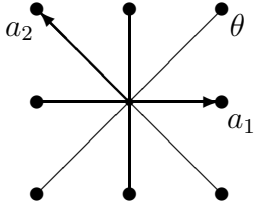
\begin{figure}
\centering
\begin{picture}(90,90)
\put(10,10){\line(1,1){70}}
\put(10,45){\line(1,0){70}}
\put(10,80){\line(1,-1){70}}
\put(45,10){\line(0,1){70}}
\put(45,45){\circle*{3}}
\put(10,10){\circle*{5}}
\put(10,45){\circle*{5}}
\put(45,10){\circle*{5}}
\put(80,10){\circle*{5}}
\put(45,80){\circle*{5}}
\put(83,70){$\theta$}
\put(80,45){\circle*{5}}
\put(82,35){$a_1$}
\put(80,80){\circle*{5}}
\put(10,80){\circle*{5}}
\put(-2,70){$a_2$}
\thicklines
\put(45,45){\vector(1,0){34}}
\put(45,45){\vector(-1,1){34}}
\end{picture}
\caption{\label{c2roots}\sl The root system of type $C_2$, with simple roots
  labeled and indicated by arrows. The lattice they generate is a scaled
  version of the ring of
  Gaussian integers. There are both long and short roots, and the same system
  gives type $B_2$.}
\end{figure}

The root system of type $C_2$ is shown in fig.~\ref{c2roots}. It is the 
same as the system of type $B_2$ but the over-extension constructed from 
this root system is $C_2^{++}$ with Dynkin diagram shown in fig.~\ref{c2pp}. 
The dual diagram corresponds to the extension of the twisted affine Lie
algebra $D_2^{(2)}$ (not $B_2^{++}$!) and will be discussed in Appendix A.

The $C_2$ root system is not simply laced, having simple roots whose 
squared lengths 
are in the ratio $2$ to $1$:
\be\label{typeC2}
\ve_1 = a_1\sqrt{2} = 1 \;\; , \; \ve_2 = a_2 =  \frac{-1 + i}{\sqrt{2}}
  \;\; , \;\; \theta = - \bar{\ve}_2 =\frac{1+i}{\sqrt{2}}\;\; .
\ee
We obtain the hyperbolic $C_2^{++}$ Dynkin diagram since our simple roots
satisfy  
\be
(\a_{-1},\a_{-1}) = (\a_0,\a_0) = (\a_2,\a_2) = 2\;\; , \; (\a_1,\a_1) = 1,
\ee 
\be
(\a_{-1},\a_0) = -1\;\; , \; (\a_0,\a_1) = -1\;\; , \; (\a_1,\a_2) = -1
\ee
and all other inner products are zero. From (\ref{typeC2}) we see 
that $\theta$ is a primitive eighth root of unity with $\theta^2= -
\bar\theta^2 = i $. The hyperbolic simple roots of $C_2^{++}$ from
(\ref{typeC2}) are
\begin{align}
\a_{-1}&=\lp\begin{array}{cc}1&0\\0&-1\end{array}\rp\,,&
\a_0&=\lp\begin{array}{cc}-1&\frac{-1-i}{\sqrt{2}}\\
   \frac{-1+i}{\sqrt{2}}&0\end{array}\rp\,,&\nn\\  
\a_{1}&=\lp\begin{array}{cc}0&\frac{1}{\sqrt{2}}\\
   \frac{1}{\sqrt{2}}&0\end{array}\rp\,,&
\a_{2}&=\lp\begin{array}{cc}0&\frac{-1+i}{\sqrt{2}}\\
         \frac{-1-i}{\sqrt{2}}&0\end{array}\rp\,.&
\end{align}
The long root $\a_2$ has $a_2$ as a unit whereas $a_1$, entering the short
root $\a_1$, is not a unit.   

To determine $W^+(C_2^{++})$ it proves convenient to bring the matrices 
$\{S_0,S_i\}$ to another form by means of a similarity transformation
\be\label{simtrm}
\tilde{S} = U S U^{-1}\;\; , \quad    U=
\lp\begin{array}{cc}\theta^{1/2}&0\\0&\theta^{-1/2}\end{array}\rp
\ee
which gives
\be\label{tS}
\tilde{S}_0 =  
\lp\begin{array}{cc}0&\theta^2\\ -\bar\theta^2&1\end{array}\rp\,,\qquad
\tilde{S}_1 = 
\lp\begin{array}{cc}0&-\theta\\\bth&0\end{array}\rp\,,\qquad
\tilde{S}_2 = \lp\begin{array}{cc}0&1\\-1&0\end{array}\rp\,.
\ee
Here, we have used $\theta=-\bar{\ve}_2$.
{}From these we can build the matrices
\begin{align}\label{psl}
\ca = \tS_2 
   &=  \lp\begin{array}{cc}0&1\\-1&0\end{array}\rp\; , &
\cb = \tS_1\tS_2
  &= \lp\begin{array}{cc}\theta&0\\0&\bar\theta\end{array}\rp\; , \nn\\
\cc = \tS_1\tS_2\tS_1 
   &=  \lp\begin{array}{cc}0&-\theta^2\\ \bar\theta^2&0\end{array}\rp\; , &
\cd = \tS_1\tS_0\tS_1\tS_2
   &=  \lp\begin{array}{cc}1&-1\\0&1\end{array}\rp\,,\nn\\
\ce = \tS_1\tS_2\tS_1\tS_0
   &=  \lp\begin{array}{cc}1&- \theta^2\\0&1\end{array}\rp\; .& \quad& 
\end{align}
The group generated by these matrices is isomorphic to the even part 
of the hyperbolic Weyl group. Hence, 
these matrices contain inversions and rotations (generated by $\ca$, $\cb$ 
and $\cc$), and translations along some lattice directions (generated by 
$\cd$ and $\ce$). The similarity transformation (\ref{simtrm}) is useful 
for explicitly exhibiting the correct lattice translations along 
two independent basis vectors $1$ and $\theta^2$ of the chosen 
integers via the matrices $\cd$ and $\ce$, respectively. We note that 
for all $\theta$ the relation $C\cdot A= B^2$ is valid, showing that 
the group generated by $A$, $B$, $C$, $D$ and $E$ is an index $2$ 
extension of the group generated by $A$, $C$, $D$ and $E$. 

With this we can easily recover the link with the so-called Klein--Fricke 
group which was first noticed in \cite{FeFr83}. Because $\theta^2=i$ 
the matrices $\cc$ and $\ce$ of (\ref{psl}) become
\be\label{psl2}
\cc =  \lp\begin{array}{cc}0&-i\\-i&0\end{array}\rp\quad , \quad
\ce = \lp\begin{array}{cc}1&-i\\0&1\end{array}\rp\; .
\ee
Together with $\ca$ and $\cd$, these matrices generate the Picard group
$PSL_2(\cG)$ (where $\cG \equiv \ints(i)$ are the Gaussian integers),
see e.g.~\cite{LaRo04}.  
We thus recover the result of \cite{FeFr83} (where, 
however, the explicit form of the embedding was not given).

\begin{prop} The even Weyl group $W^+(C_2^{++})$ is an index 2 extension of
$PSL_2(\cG)$, that is, 
\be
W^+(C_2^{++}) \cong PSL_2(\cG)\rtimes 2 \equiv PSL_2(\ints(i))\rtimes 2 \,.
\ee 
\end{prop}
Semi-directness follows since conjugation of $PSL_2(\cG)$ by the matrix $B$ is
an automorphism. 

\begin{figure}
\centering
\begin{picture}(110,20)
\thicklines
\multiput(10,10)(30,0){4}{\circle*{8}}
\put(14,10){\line(1,0){22}}
\put(40,14){\line(1,0){30}}
\put(40,6){\line(1,0){30}}
\put(70,14){\line(1,0){30}}
\put(70,6){\line(1,0){30}}
\put(50,20){\line(1,-1){10}}
\put(50,0){\line(1,1){10}}
\put(90,20){\line(-1,-1){10}}
\put(90,0){\line(-1,1){10}}
\put(4,-5){$-1$}
\put(37,-5){$0$}
\put(67,-5){$1$}
\put(97,-5){$2$}
\end{picture}
\caption{\label{c2pp}\sl $C_2^{++}$ Dynkin diagram with numbering of
 nodes.}
\end{figure}
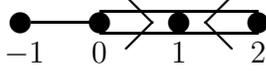

\end{subsection}

\begin{subsection}{$\Kk=\cx$, type $G_2$}

The root system of type $G_2$ is depicted in fig.~\ref{g2roots}. It is 
not simply laced, having simple roots whose squared lengths are in 
the ratio $3$ to $1$.
\be\label{typeG2}
\ve_1 = a_1 = 1\;\; , \; \ve_2=\sqrt{3} a_2 = \frac{-\sqrt{3}+i}{2} \;\; , \; 
 \theta = \frac{1+\sqrt{3}i}{2}\;.
\ee
The $G_2$ root system is thus the superposition of two $A_2$ root systems
which are scaled by a factor of $\sqrt{3}$ and rotated by $30^\circ$ degrees
relative to each other. Note that in this case $\theta\ne -\bar{\ve}_2$;
rather, we have $\ve_2^2= \bth$, and that $\theta$ is identical to the
highest root of $A_2$.

The Dynkin diagram of the hyperbolic algebra $G_2^{++}$ is shown in 
fig.~\ref{g2pp}. We will get this diagram if our simple roots satisfy 
\be
(\a_{-1},\a_{-1}) = (\a_0,\a_0) = (\a_1,\a_1) = 2\;\; , 
\; (\a_2,\a_2) = \frac{2}{3},
\ee 
\be
(\a_{-1},\a_0) = -1\;\; , \; (\a_0,\a_1) = -1\;\; , \; (\a_1,\a_2) = -1
\ee
and all others zero. 
The hyperbolic simple roots from (\ref{typeG2}) of $G_2^{++}$ are
\begin{align}
\a_{-1}&=\lp\begin{array}{cc}1&0\\0&-1\end{array}\rp\,,&
\a_0&=\lp\begin{array}{cc}-1&\frac{-1-\sqrt{3}i}{2}\\
    \frac{-1+\sqrt{3}i}{2}&0\end{array}\rp\,,&\nn\\
\a_{1}&=\lp\begin{array}{cc}0&1\\1&0\end{array}\rp\,,&
\a_{2}&=\lp\begin{array}{cc}0&\frac{-1+\sqrt{3}i}{2}\\
     \frac{-1-\sqrt{3}i}{2}&0\end{array}\rp\,.&
\end{align}
(NOTE: This formula for $\a_2$ had a typographical error in the published 
version which is corrected above.) 
We label $\a_2$ as a short root, so $a_2$ is not a unit. 

The finite Weyl group $W(A_2)$ is a subgroup of index 2 in $W(G_2)$.
The same is true for their hyperbolic extensions: we have

\begin{prop} The even Weyl group $W^+(G_2^{++})$ is an index 2 extension
of $W(A_2^{++})$, that is,  
\be
W^+(G_2^{++}) = W^+(A_2^{++})\, \rtimes\,2 = PSL_2(\cE)\,\rtimes\,2 \,.
\ee
\end{prop}

\noindent{\bf{Proof:}} Because $\theta \neq -\bve_2$ we must proceed 
slightly differently than before. First we notice that
\be
S_1 S_2 = \lp\begin{array}{cc}- \bve_2&0\\0& - \ve_2\end{array}\rp \;\;
\Rightarrow \quad  \big(S_1 S_2\big)^2 = 
\lp\begin{array}{cc} \theta &0\\0& \bth\end{array}\rp \;
\ee
Since $\theta^6=1$ (as for $A_2$), the set of matrices $S_0,S_1$ and
$(S_1S_2)^2 S_1$ coincides with the set (\ref{SC}) for $A_2$,
hence these matrices generate the group $PSL_2(\cE)$, again by
Proposition~\ref{sl2prop}. To get the full even Weyl group we must adjoin 
the matrix $S_2$ obeying $S_2^2 = - \id$, generating a $\Ztwo$. 
The semi-directness of the product will follow from the action $S_2 S
S_2^{-1}\cong S_2 S S_2$ of the extending matrix $S_2$ on $S\in
PSL_2(\cE)$, if the resulting matrix is in $PSL_2(\cE)$ so that this
action gives an automorphism of $PSL_2(\cE)$. By expanding the
product, commutativity and the fact that $\ve_2^2\in \cE$ the result
follows.$\square$

\begin{figure}
\centering
\begin{picture}(100,140)
\put(115,48){$a_1$}
\put(-14,75){$a_2$} 
\put(75,125){$\theta$}

\put(40,60){\line(0,1){40}} 
\put(40,60){\line(0,-1){40}} 
\put(40,60){\line(2,1){40}} 
\put(40,60){\line(-2,1){40}} 
\put(40,60){\line(2,-1){40}} 
\put(40,60){\line(-2,-1){40}} 
\put(80,40){\circle*{5}} 
\put(80,80){\circle*{5}} 
\put(40,20){\circle*{5}} 
\put(40,100){\circle*{5}} 
\put(0,40){\circle*{5}} 
\put(0,80){\circle*{5}} 

\put(40,60){\line(1,0){80}}
\put(40,60){\line(-1,0){80}}
\put(40,60){\line(2,3){40}}
\put(40,60){\line(-2,3){40}}
\put(40,60){\line(2,-3){40}}
\put(40,60){\line(-2,-3){40}}
\put(-40,60){\circle*{5}} 
\put(120,60){\circle*{5}} 
\put(0,120){\circle*{5}} 
\put(0,0){\circle*{5}} 
\put(80,120){\circle*{5}} 
\put(80,0){\circle*{5}} 
\thicklines
\put(40,60){\vector(1,0){80}}
\put(40,60){\vector(-2,1){40}}
\end{picture}
\caption{\label{g2roots}\sl The root system of type $G_2$. The long and 
short roots are labeled and indicated by arrows.} 
\end{figure}
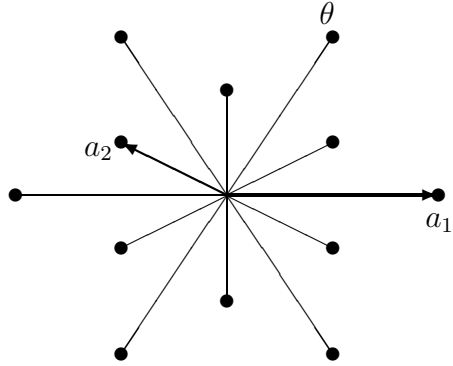

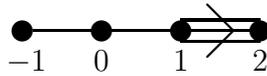
\begin{figure}
\centering
\begin{picture}(110,20)
\thicklines
\multiput(10,10)(30,0){4}{\circle*{8}}
\put(14,10){\line(1,0){22}}
\put(40,10){\line(1,0){30}}
\put(70,14){\line(1,0){30}}
\put(70,10){\line(1,0){30}}
\put(70,6){\line(1,0){30}}
\put(80,20){\line(1,-1){10}}
\put(80,0){\line(1,1){10}}
\put(4,-5){$-1$}
\put(37,-5){$0$}
\put(67,-5){$1$}
\put(97,-5){$2$}
\end{picture}
\caption{\label{g2pp}\sl $G_2^{++}$ Dynkin diagram with numbering of nodes.}
\end{figure}

\end{subsection}

\end{section}

\begin{section}{Quaternions $\Kk = \quat$}
\label{quatsec}

Within the four-dimensional quaternion algebra one can find the root systems
of types $A_4$, $B_4$, $C_4$, $D_4$ and $F_4$. The associated hyperbolic
rank $6$ Kac--Moody algebras are now $A_4^{++}, B_4^{++}, C_4^{++}, D_4^{++}$,
and $F_4^{++}$, and we will give explicit descriptions of their even
Weyl groups in terms of matrix groups below.\footnote{A description of the 
finite root systems of types $B_4$, $D_4$ and $F_4$, and their Weyl groups
in terms of quaternions has been given in \cite{KoKoAl03}.}
The quaternionic case is more subtle than the commutative cases because 
the criterion for selecting the matrix group to which the even Weyl group 
belongs cannot be so easily done via a determinant. Before turning to the 
issue of how to define determinants and matrix groups with quaternionic 
entries we first discuss different types of quaternionic integers and numbers
required for exhibiting the root systems. 

\begin{subsection}{Hurwitz and Lipschitz integers}

The standard basis for $\quat$, $\{1, i, j, k\}$, has the famous products 
\be
i^2 = j^2 = k^2 = -1,\ \  ij=-ji=k,\ \  jk = -kj = i,\ \  ki = -ik = j,
\ee
so an obvious subring of integers is formed by the {\em Lipschitz integers} 
\be\label{Lnum}
\Lq = \left\{n_0 + n_1 i + n_2 j + n_3 k \;\vline\; n_0, n_1, n_2, n_3\in
 \ints\right\}\,. 
\ee
They form an order (a subring of the rational quaternions, $\quat_\rats$, 
finitely generated as a $\ints$-module, containing a $\rats$-basis of
$\quat_\rats$) but are not a maximal order in $\quat_\rats$ since they
are contained in the ring of {\em Hurwitz integers}
\begin{equation}\label{Hnum}
\Hq = \left\{n_0 + n_1 i + n_2 j + n_3 k\; \vline\; n_0, n_1, n_2, n_3\in
 \ints\;\text{or}\; 
n_0, n_1, n_2, n_3\in \ints+\frac12 \right\}\,,
\end{equation}
which constitute a maximal order in $\quat_\rats$. We note that these
two rings, $\Lq$ and $\Hq$, are generated by the Lipschitz and Hurwitz units 
given below in (\ref{LipUnits}) and (\ref{HurUnits}), respectively, and as 
$\ints$-modules, they are discrete lattices. 
Also note that the non-commutative ring $\Hq$
of Hurwitz integers satisfies the division with small remainder property 
required for the Euclidean algorithm (cf.~Proposition~\ref{sl2prop}), whereas 
the Lipschitz integers $\Lq$ do not \cite{CoSm03}. 

For the determination of some of the even Weyl groups we rely on the following
definition and lemma.

Define $\Cq$ to be the two-sided ideal in the ring $\Hq$ generated by
the commutators 
\be
[a,b] = ab - ba \qquad\hbox{ for all }a,b\in\Hq\,.
\ee
To understand this we use the integral basis for $\Hq$ used in \cite{Kr85},
\be
e_0 = \frac12(1+i+j+k), \qquad i, \qquad j,\qquad k,
\ee
which has the nice property that an integral combination 
$m_0 e_0 + m_1 i + m_2 j + m_3 k$ is in $\Lq$ when $m_0\in 2\ints$. 
It is easy to check the commutators
\begin{align} \label{Hcomm}
[e_0,i] &= j-k, &  [e_0,j] &= k-i,& [e_0,k] &= i-j, \\ \nn 
[i,j] &= 2k, &  [j,k] &= 2i, & [k,i] &= 2j,
\end{align}
which are all purely imaginary Lipschitz integers. $[\Hq,\Hq]$ is all integral 
linear combinations of the commutators above, but that is not even a
subring of  $\Hq$, as we can see, for example, from the fact that
$(2i)(j-i) - 2k = 2$ is not a commutator. The ideal $\Cq = \Hq
[\Hq,\Hq] \Hq$ consists of finite sums of the form $a [b,c] d$ for any
$a,b,c,d\in\Hq$, and it is enough to compute these for 
$a$ and $d$ from the above integral basis and $[b,c]$ from (\ref{Hcomm}). 
The only products we need to know are
\begin{align}
e_0 i &= -e_0+i+j, & e_0 j &= -e_0+j+k, & e_0 k &= -e_0+i+k, \nn\\
i e_0 &= -e_0+i+k, & j e_0 &= -e_0+i+j, & k e_0 &= -e_0+j+k, \nn\\
 & & e_0 e_0 &= -e_0+i+j+k. & & 
\end{align}
Using these, it is clear that $a [b,c] d$ is in $\Lq$, and consists of all
integral linear combinations of the elements
\be
\pm 1\pm i,\ \ \pm 1\pm j,\ \ \pm 1\pm k,\ \ \pm i\pm j,\ \ \pm i\pm
k,\ \ \pm j\pm k,  
\ee
which is clearly an index $2$ integral lattice in $\Lq$. Since $\Lq$
is an index $2$ sublattice in $\Hq$, we have that $\Hq/\Cq$ is a ring
of four elements, and we can take as coset representatives $\{0, -e_0,
(-e_0)^2 = e_0-1, (-e_0)^3 = 1\}$ which form the field $\Ff_4$ of
order $4$ whose nonzero elements form the cyclic group of order $3$.  
Also note that there are no units in $\Cq$ since its nonzero elements
have minimal length $2$. 

\begin{lemma}\label{HL} 
Let $a_1,\ldots, a_n\in\Hq$ be any $n$ Hurwitz numbers. Then the product 
$a_1\cdots a_n$ is commutative modulo $\Cq$.
\end{lemma}

\noindent
{\bf{Proof:}} Since the quotient ring $\Hq/\Cq$ is a field, where the
product is  commutative, the projection of the product $a_1\cdots a_n$
is equal to the projection of $a_{\sigma(1)}\cdots a_{\sigma(n)}$ for
any permutation $\sigma$ of $\{1,\ldots,n\}$. $\square$

\end{subsection}

\begin{subsection}{Quaternionic units and rings}
\label{quatunits}

For the root lattices of $A_4$, $B_4$, $C_4$ and $F_4$ (but not $D_4$) 
we also need quaternionic units which are neither Lipschitz nor Hurwitz 
numbers. These other units parametrize finite subgroups of $SU(2)$, and 
are, respectively, related to the octahedral (for $B_4$, $C_4$ and
$F_4$) and icosahedral (for $A_4$)
groups, as explained e.g. in \cite{CoSm03}. Using diagonal and off-diagonal 
$(2\times 2)$ matrices of such units (or alternatively {\em pairs} of units
\cite{CoSm03}, cf. footnote~\ref{csnot} in Section~\ref{faff}) we can then
reconstruct the Weyl groups of all the finite simple rank 4 algebras,
as we shall explain below. Besides the eight Lipschitz units 
\be\label{LipUnits}
\Units_\Lq = \{ \pm 1,\pm i, \pm j, \pm k \}
\ee 
which form the quaternionic group, often denoted by $Q_8$, we have the
$24$ Hurwitz units 
\be\label{HurUnits}
\Units_\Hq = \left\{ \pm 1 \,,\,\pm i\,,\, \pm j\,,\, \pm k \,,\, 
   \frac12(\pm 1 \pm i \pm j \pm k) \right\}
\ee
which form a subgroup in the unit quaternions. According to
\cite{CoSl88}, p.~55, $\Units_\Hq \cong 2\cdot \Ag_4$ is an index two
extension of the alternating group on four letters. It can also be
seen that $\Units_\Hq \cong \Units_\Lq \rtimes \ints_3$ is a
semi-direct product of the group of Lipschitz units with a cyclic group
of order $3$ related to triality and explained below in Section~\ref{d4sec}.
Any quaternion $z=n_0 + n_1 i + n_2 j + n_3 k$ is a root of the real polynomial 
\be
p_z(t) = (t - z)(t - {\bar{z}}) = t^2 - 2n_0 t + (n_0^2+n_1^2+n_2^2+n_3^2)
\ee
whose coefficients will be integers if $z$ is in $\Lq$ or $\Hq$, and
in that case, unless $z\in\ints$, this is the minimal polynomial
satisfied by $z\in\Hq$.  

We will also need the following set of $24$ {\em octahedral} units 
\be
\left\{ \frac{\pm 1 \pm i}{\sqrt{2}} \,,\,
\frac{\pm 1 \pm j}{\sqrt{2}}\,,\,\frac{\pm1 \pm k}{\sqrt{2}}\,,\,
\frac{\pm i \pm j}{\sqrt{2}}\,,\, \frac{\pm i \pm k}{\sqrt{2}}\,,\,
\frac{\pm j \pm k}{\sqrt{2}}  \right\}
\ee
which do not form a group by themselves, but the product of any two of
them is a Hurwitz unit, and the product of any Hurwitz unit and an
octahedral unit is an octahedral unit, so their union forms a group of
order $48$ which we call the octahedral subgroup, $\Units_\cR$. We can
write it as the disjoint union of two cosets of its normal subgroup,
$\Units_\Hq$,  
\be\label{ounits}
\Units_\cR \equiv \Units_\Hq \cup i_{\text{O}} \cdot\Units_\Hq \qquad\qquad
\mbox{where}  \quad i_{\text{O}}\equiv \frac{j-k}{\sqrt{2}}\,.
\ee
These are, in fact, all the units in the ring, $\cR$, generated by
integral linear combinations of units in $\Units_\cR$. Contrary to the
Lipschitz and Hurwitz numbers, the ring $\cR$ is not discrete but {\em
  dense} in $\quat$ (using the usual topology). Nevertheless $\cR$ is
an order of $\quat_\Ff$ (with $\Ff=\rats(\sqrt{2})$) as defined in
Section~\ref{intsec}, since it is finitely generated by its $48$ units
over $\ints$. It is not hard to show that it is generated over $\Ss =
\ints[\sqrt{2}]$ just by the four elements, $\{a,b,ab,ba\}$, where $a =
\frac{1+i}{\sqrt{2}}$ and $b = \frac{1+j}{\sqrt{2}}$. Note also that
the integral span of the coset of `purely octahedral numbers'
$i_{\text{O}}\cdot\Hq$ constitutes a lattice in $\quat$ which can be
regarded as a `rotated' version of the lattice $\Hq$ (but which,
unlike $\Hq$, is not closed under multiplication).\footnote{In fact,
  we will only encounter matrices $S$ whose entries belong either 
  to $\Hq$ (and have minimal polynomial of degree two) or to $i_{\text{O}}
  \Hq$ (and have minimal polynomial of degree four). This will be distributed
  over the matrix in such a way that this structure is preserved under
  multiplication.} 

Finally, we have $96$ {\em icosian} units \cite{CoSl88}, p.~207,
defined in terms of  
\be
\tau=\frac12 (1+\sqrt{5}) \quad\text{and}\quad
\si=\frac12 (1-\sqrt{5})\,,
\ee
to be all elements obtained from the following eight basic elements
\be
\frac12 \big( \pm i \pm \si j \pm \tau k \big)
\ee
by the $12$ {\em even} permutations of the Lipschitz units $(1,i,j,k)$. 
Together with the $24$ Hurwitz units they form the icosian subgroup
$\Units_\cI$, which are the $120$ units in the ring $\cI$ defined to
be their integral span. According to \cite{CoSl88}, the structure of
the icosian group $\Units_\cI \cong 2\cdot \Ag_5$ is an index two
extension of the alternating group on five letters.  
$\cI$ is an order of $\quat_\Ff$ where $\Ff = {\rats(\sqrt{5})}$. As
in the order $\cR$, the ring $\cI$ is dense in $\quat$, but it is
finitely generated over the integers by its $120$ units. In fact, it
is generated over $\Ss = \ints(\tau)$ by the four elements 
\be 
i\ , \quad j\ , \quad \omega = \frac12 (-1+i+j+k)\ , \quad 
i_\cI = \frac12 (i + \si j + \tau k)\,. 
\ee
It is straightforward to check that all their products 
can be expressed as linear combinations of those four elements with
coefficients from $\Ss$, for example:
\begin{align}
i^2 &= -1 = -(\si+1)i-(\si+2)j+2\omega+2\si i_\cI, \\ \nn
ji_\cI &= -2\si i-(2\si+1)j+\si\omega+(2\si+1)i_\cI, & 
   \omega i_\cI &= \omega-i-k \\ \nn 
i_\cI \omega &= -i-(\si+1)j+\omega+\si i_\cI, 
   & i_\cI j &= -i-\si j+\si\omega+i_\cI\,. 
\end{align}
These can be used to compute the commutators of pairs of generators, e.g.
\begin{align}\label{Icom}
[i,j] &= 2k, & [i,i_\cI] &= -\tau j+\si k, 
\end{align}
If we try to imitate what we did for the Hurwitz integers, and let
$\Cq_\cI = \cI[\cI,\cI]\cI$ be the ideal of $\cI$ generated by the
commutators $[\cI,\cI]$, we find that $\Cq_\cI = \cI$ is the whole
icosian ring. This follows since the two elements in (\ref{Icom}) have
norms $2$ and $3$, respectively, so that the difference of their norms
is $1$ and hence $1\in \Cq_\cI$. This means that the quotient ring,
$\cI/\Cq_\cI$ is trivial, not a field as was $\Hq/\Cq$, and we do not
have an analog of Lemma~\ref{HL} for $\cI$. 

Even though most rank four root systems contain units from one of the 
dense rings $\cR$ or $\cI$ we will explain how one can nevertheless 
construct {\em discrete matrix groups} from them which are related to 
the even Weyl groups.

\end{subsection}

\begin{subsection}{Quaternionic determinants}

Useful introductory references for quaternionic determinants are 
\cite{Dy72,As96}. When $\Kk = \quat$, the matrices in (\ref{evenmat}) are 
$(2\times2)$ matrices over the quaternions, that is, matrices of the form
\be\label{psl2hmat}
S =  \lp\begin{array}{cc} a& b\\c& d \end{array}\rp \qquad 
\mbox{with $a,b,c,d \in \quat$} \,.
\ee
To describe the hyperbolic Weyl groups, we will only need to take these 
entries from the rings introduced above, but the following observations apply 
to all quaternionic matrices. If one tries to define the determinant of the 
quaternionic matrix $S$ to be $ad-bc$, the non-commutativity of $\quat$  
spoils the usual multiplication law for determinants (see \cite{As96}). 
A quantity which does have good properties for $\Kk=\quat$ is obtained by 
observing that the determinant is well defined for {\em Hermitian}
quaternionic matrices (for which $a$ and $d$ are real, and
$b=\bar{c}$). Given any quaternionic matrix $S$, we therefore
associate with it the Hermitian matrix $SS^\dagger$ which has a
well-defined determinant given by  
$\det (S S^\dagger)=\det(S^\dagger S)$.\footnote{These expressions also 
make sense for $\Kk =\oct$ but in general are no longer equal.}
The multiplication law is then obeyed provided that
\be
\det \big( (S_1 S_2)\cdot (S_2^\dagger S_1^\dagger) \big) 
= \det \big( S_1 S_1^\dagger\big) \, \det \big(S_2 S_2^\dagger \big)\,.
\ee
This relation (which is obviously satisfied for $\Kk=\reals,\cx$) can  
be verified for $\Kk = \quat$ by direct calculation (however, it will fail 
for $\Kk=\oct$ because of non-associativity). Any quaternionic matrix $S$ 
satisfying $\det(SS^\dagger)\neq 0$ has left and right inverses given by 
\be\label{lrinverses}
\big( S^\dagger S \big)^{-1} S^\dagger = 
S^\dagger \big( S S^\dagger \big)^{-1}\,.
\ee
Here, $(S S^\dagger)^{-1}$ denotes the well-defined inverse of a Hermitian
matrix. That the two expressions coincide can be shown by explicit 
computation.\footnote{For $\Kk=\oct$, the two expressions in
  (\ref{lrinverses}) are  
  generally different because of non-associativity.} 
Hence, the condition
$\det \big( S S^\dagger \big) = 1$ for $(2 \times 2)$ matrices with
quaternionic entries indeed defines a continuous group which we denote by
\be\label{PSLH}
SL_2(\quat) := \left\{ S =  
 \lp\begin{array}{cc} a & b\\ c & d \end{array}\rp 
 \;\; \vline\;\; a,b,c,d \in \quat \; ,\;
  \det \big( SS^\dagger \big) = 1 \right\}\,.
\ee 
Writing $\id$ for the $(2\times 2)$ identity matrix, the subgroup
$\{\pm \id\}$ is normal and we define $PSL_2(\quat) = SL_2(\quat)/\{\pm \id\}$
to be the quotient group. We write elements in the quotient as matrices $S$
but we identify $S$ with $-S$.\footnote{In fact, if $S\in SL_2(\quat)$ then 
  so is $\ve S$ for any unit $\ve$. Because of non-commutativity, only 
  $\ve=\pm1$ are in the center and quotiented out.}
This group has real dimension $15 = 4\cdot 4 -1$, and is known to be 
isomorphic to the Lorentz group in six dimensions \cite{Su84,Ba02}, i.e. 
\be\label{pso15}
PSL_2(\quat) \cong PSO(1,5;\reals)\,.
\ee
We will soon see that the even Weyl groups of all rank six Kac--Moody algebras 
under consideration are discrete subgroups of $PSL_2(\quat)$.

The easiest of these discrete groups to describe is $PSL_2(\Hq)$, 
which is obtained from $PSL_2(\quat)$ by restricting all matrix entries 
to be Hurwitz integers. We will also need its subgroup 
\be\label{gq}
\PGL := \Big\{S = \lp\begin{array}{cc} a & b\\ c & d \end{array}\rp
\in PSL_2(\Hq)\;\vline\; ad-bc \equiv 1 (\text{mod}\,\Cq) \Big\}
\ee 
which can be understood in the following way. The composition of the group 
homomorphism
\be
PSL_2(\Hq) \to PSL_2(\Hq/\Cq)
\ee
with the usual determinant of a $(2\times 2)$ matrix over the field 
$\Hq/\Cq \cong \Ff_4$, yields a group homomorphism
\be
{\rm Det}: PSL_2(\Hq) \to \Ff_4^*
\ee
onto $\Ff_4^* \cong \ints_3$ the cyclic group of order 
$3$.\footnote{This `determinant' is different from the 
  Dieudonn\'e{} determinant \cite{As96}.}
Since this composition is a group homomorphism, we have 
${\rm Det}\, (S_1 S_2) = {\rm Det} \, S_1 \cdot {\rm Det} \, S_2$. 
Its kernel is a normal subgroup of index $3$ in $PSL_2(\Hq)$, giving
an alternative form of (\ref{gq}) as
\be\label{gq2}
\PGL = \big\{ S\in PSL_2(\Hq) \,\,\vline\,\, \rmD\, S = 1 \big\}.
\ee 
Since $\Lq = \Cq \cup (1+\Cq)$ is the coset decomposition of $\Lq$, the condition
$ad-bc \equiv 1 (\text{mod}\,\Cq)$ says that $ad-bc\in (1+\Cq)\subset\Lq$
and Lemma~\ref{HL} says that the order of the products in the expression  
$ad - bc$ does not matter. Since $ad-bc\notin\Cq$, 
the condition for $S\in\PGL$ is just $ad-bc\in\Lq$. We have proved the 
following. 

\begin{lemma}\label{d4lemma}
$\PGL$ is index $3$ in $PSL_2(\Hq)$, and
\be
PSL_2(\Hq)/\PGL = {\mathfrak A}_3
\ee
where ${\mathfrak A}_3\cong \ints_3$ is the alternating group on three
letters. 
\end{lemma}

The modular group $\PGL$ contains the modular group $PSL_2(\Lq)$, but 
is strictly larger, so that we arrive at the following chain of 
subgroup relations
\be
PSL_2(\Lq)\subset\PGL\subset PSL_2(\Hq)\,.
\ee

The map Det extends to the ring of all $(2\times 2)$ Hurwitz matrices, first
reducing entries modulo $\Cq$, and then taking the usual determinant. It is 
still a multiplicative map so that the invertible matrices have non-zero Det.
But it is possible for a `non-invertible' matrix in that ring to
also have a non-zero Det. For example, the diagonal matrix
$\text{diag}(1,3)$ does not have an inverse over the Hurwitz integers,
but since $2\in\Cq$, it reduces to the identity matrix modulo $\Cq$,
whose usual determinant is $1+\Cq\in\Ff_4^*$.

\begin{prop}\label{psl0gen}
The modular group $\PGL$ is generated by the matrices 
\be\label{gl2hgen}
\lp\begin{array}{cc}1&1\\0&1\end{array}\rp\,,\quad
\lp\begin{array}{cc}0&1\\1&0\end{array}\rp\,,\quad
\lp\begin{array}{cc}a &0\\0&b\end{array}\rp\,,
\ee
where $a,b\in\Units_\Hq$ are Hurwitz units and $ab\equiv 1 (\text{mod}\ \Cq)$. 
\end{prop}

\noindent{\bf{Proof:}}
This claim can be proved by adapting Theorem~2.2 on page~16 of \cite{Kr85}, 
and arguments very similar to those of Proposition~\ref{sl2prop}. According 
to Krieg's theorem, the set of all uni-modular (= invertible) $(2 \times 2)$
Hurwitz quaternionic matrices, here denoted by $SL_2(\Hq)\subset SL_2(\quat)$ in
accordance with definition (\ref{PSLH}) (but designated as $GL_2(\Hq)$ in 
\cite{Kr85}) is generated by taking products of the matrices in (\ref{gl2hgen})
without any restriction on the product $ab$, so that $a,b$ run through
all Hurwitz units. This is also true for the quotient $PSL_2(\Hq)$. So we can
write any element $S\in PSL_2(\Hq)$ as 
\be\label{kr}
S=G_1\cdots G_n
\ee
where each $G_i$, $i=1,\ldots,n$, is one of the three types of
generating matrices in (\ref{gl2hgen}), but with no restriction on the
product $ab$ for the third type in accordance with Krieg's theorem. We
will show that we can rewrite $S$ in the form
\be\label{krnormal}
S = \lp\begin{array}{cc}c&0\\0&1\end{array}\rp G'_1\cdots G'_{n'}\,,
\ee
where $c\in\Hq$ and all $G'_i$ matrices are from (\ref{gl2hgen}) with
the third type satisfying the $ab\equiv 1 (\text{mod}\ \Cq)$ condition. 
In order to arrive at this new presentation we show
explicitly how to convert the presentation (\ref{kr}) to one of the form 
(\ref{krnormal}), and how the element $c$ represents $\rmD(S)$.

Examine the matrices in the expression $G_1\cdots G_n$ from the right
to the left. We don't need to change a matrix of the first or second
kind from (\ref{gl2hgen}), nor do we change one of the third kind if
it satisfies the $ab\equiv 1 (\text{mod}\ \Cq)$ condition. Let $i_0$
be the largest index for which $G_{i_0}=\text{diag}(a_{i_0},b_{i_0})$
with $a_{i_0}b_{i_0}\not\equiv 1 (\text{mod}\ \Cq)$. Then we can
factorize it as  
\be
\lp\begin{array}{cc}a_{i_0}&0\\0&b_{i_0}\end{array}\rp
= \lp\begin{array}{cc}\bar{\ve}&0\\0&1\end{array}\rp
\lp\begin{array}{cc}\ve a_{i_0}&0\\0&b_{i_0}\end{array}\rp 
\ee
where the unit $\ve\in\Units_\Hq \setminus \Units_\Lq$ is chosen such that 
$\ve a_{i_0}b_{i_0} = 1 \in \Lq$. We will show below how the new 
matrix $\text{diag}(\bar{\ve},1)$ can be moved to the left by passing 
through the other generators $G_i$ with $i<i_0$, and leaving only
acceptable generators on its right.

We first note that from (\ref{gl2hgen}) (with the restriction on 
the product) we can generate all translation matrices 
\be\label{trq}
\lp\begin{array}{cc}1&a\\0&1\end{array}\rp\,,\quad \text{for}\, a\in \Hq\,.
\ee
For $a\in\{\pm i,\pm j,\pm k\}$ we have
\be
\lp\begin{array}{cc}a&0\\0&1\end{array}\rp
\lp\begin{array}{cc}1&1\\0&1\end{array}\rp
\lp\begin{array}{cc}-a&0\\0&1\end{array}\rp
=\lp\begin{array}{cc}1&a\\0&1\end{array}\rp
\ee
and $a\equiv 1\ ({\text mod}\ \Cq)$ so each matrix on the left side is
from (\ref{gl2hgen}). Products of these provide 
all translation matrices (\ref{trq}) with $a\in\Lq$. To extend
this to all $a\in\Hq$ it is sufficient to find one translation
matrix with $a\in\Hq\setminus\Lq$. The computation in (\ref{trmat})
where $\theta = \frac12\left(1-i-j-k\right)$ satisfies
$\theta^2=-\bar\theta = \frac12\left(-1-i-j-k\right)$, yields this, so
that we indeed obtain all matrices (\ref{trq}) from
(\ref{gl2hgen}). From this and 
\be\label{move1}
\lp\begin{array}{cc}1&1\\0&1\end{array}\rp
\lp\begin{array}{cc}\bar{\ve}&0\\0&1\end{array}\rp
=\lp\begin{array}{cc}\bar{\ve}&0\\0&1\end{array}\rp
\lp\begin{array}{cc}1& {\ve}\\0&1\end{array}\rp
\ee
we see that we can pass the matrix $\text{diag}(\bar{\ve},1)$ across a 
matrix $G_i$ of the first kind in (\ref{gl2hgen}) by replacing $G_i$ by 
a translation matrix which is a product of the allowed generators. 

Similarly we have
\be\label{move2}
\lp\begin{array}{cc}0&1\\1&0\end{array}\rp
\lp\begin{array}{cc}\bar{\ve}&0\\0&1\end{array}\rp
=\lp\begin{array}{cc}\bar{\ve}&0\\0&1\end{array}\rp
\lp\begin{array}{cc}\ve&0\\0&\bar{\ve}\end{array}\rp
\lp\begin{array}{cc}0&1\\1&0\end{array}\rp\,,
\ee
so that we can also pass through the second type of generating
matrices by replacing it by a product of allowed matrices. 

Finally, for any units $a,b\in\Hq$, we can also write
\be\label{move3}
\lp\begin{array}{cc}a&0\\0&b\end{array}\rp
\lp\begin{array}{cc}\bar{\ve}&0\\0&1\end{array}\rp
=\lp\begin{array}{cc}a\bar{\ve}b&0\\0&1\end{array}\rp
\lp\begin{array}{cc}\bar{b}&0\\0&b\end{array}\rp\,,
\ee
so the matrix on the far right is an allowed generator, but the 
matrix being moved to the left has been changed to another of the 
form $\text{diag}(c,1)$ with $c = a\bar{\ve}b\in\Hq$ a unit. Using 
(\ref{move1})--(\ref{move3}), for some unit $c\in\Hq$, we finally arrive 
at the expression 
\be
S=\lp\begin{array}{cc}c&0\\0&1\end{array}\rp G'_1\cdots G'_{n'}
\ee
with all $G'_i$ belonging to the list of generating matrices (\ref{gl2hgen})
and those of the third kind have $ab\equiv 1 (\text{mod}\ \Cq)$, so  all 
$\rmD(G'_i) = 1\in\Ff_4^*$. Then
$\rmD(S) = \rmD(\text{diag}(c,1)) \rmD(G'_1)\cdots \rmD(G'_{n'}) = c+\Cq \in 
\Hq/\Cq = \Ff_4$. If we assume that
$S\in\PGL$, so $\rmD(S) = 1$, then we must have $c \equiv 1 (\text{mod}\ \Cq)$
so that $\text{diag}(c,1)$ is an allowed generator from (\ref{gl2hgen}) 
and the proof is complete. $\square$\medskip

We close this section by counting the number of diagonal matrices 
$\text{diag}(a,b)$ such that $a,b\in \Hq$ are units and $ab\in\Lq$. There 
are $64=8\times 8$ such matrices when both $a,b\in\Lq$. If 
$a\in\Hq\backslash\Lq$ is one of the $16$ pure Hurwitz units there 
are $8$ choices for $b$ such that $ab\in\Lq$, which gives $128=16\times 8$ 
such matrices. Descending to the quotient by $\{\id,-\id\}$ we are left 
with $96 = \frac12\times 192$ diagonal matrices in $\PGL$. The full 
$PSL_2(\Hq)$ has $288 = \frac12 \times 24 \times 24$ diagonal matrices. 
This gives another way to see Lemma~\ref{d4lemma}.

\end{subsection}

\begin{subsection}{Even Weyl group $W^+(D_4^{++})$}
\label{d4sec}

We first study the Weyl group of the rank 6 hyperbolic algebra 
$D_4^{++}$. Among the Weyl groups associated to the rank 6 algebras,
this is the `easiest' to understand because it is the only one
whose root system can be expressed solely in terms of Hurwitz numbers. 
As we will see, the Weyl groups in all the other cases, except for $A_4^{++}$, 
are finite extensions of this one.

In the lattice $\Hq$ equipped with the inner product (\ref{bilinear}), we can 
recognize that $\Units_\Hq$ forms the $D_4$ root system, and among these we 
choose as simple roots
\begin{align}\label{D4roots}
\ve_1 &= a_1 = 1 \;\; , & 
\ve_2 &= a_2 = \frac12 \big( -1 + i - j - k \big) \; ,\\
\ve_3 &= a_3 =\frac12 \big( -1 - i + j - k \big) \; , &
\ve_4 &= a_4 = \frac12 \big( -1 - i - j + k \big)\,.\nn
\end{align}
We can choose $a_i=\ve_i$ for all $i$ because $D_4$ is simply
laced. This system  
of simple roots exhibits $\Sg_3$ `triality' symmetry as outer
automorphisms, that 
is, it is symmetric under cyclic permutation of the three imaginary
units $(i,j,k)$, 
and under the exchange of any two imaginary units,
e.g. $(j\leftrightarrow k)$.  
The former is concretely realized by the map $z\mapsto\theta
z\bar{\theta}$, whereas  
the latter corresponds to $z\mapsto i_{\text{O}} \bar{z}
\bar{\imath}_{\text{O}}$.  
Here, $\theta$ is the highest $D_4$ root
\be\label{hstd4} 
\theta = 2a_1 + a_2 + a_3 + a_4 = \frac12 \big( 1 - i - j - k \big)  
 = j \ve_2 = k \ve_3 = i \ve_4  
\ee
and $i_{\text{O}}$ is the specific octahedral unit of order four
defined in (\ref{ounits}).
The 16 units in $\Units_\Hq\setminus\Units_\Lq$ are of order three or
six, {\it viz.} 
\be\label{3}
\ve_2^3 = \ve_3^3 = \ve_4^3 = 1\,,\quad \theta^6 = 1\,.
\ee
With these choices we build the simple roots (\ref{simroots}) of the
hyperbolic algebra  $D_4^{++}$ in its root lattice $\Lambda(Q)$ where
$Q = \Hq$ is the $D_4$ root lattice. The corresponding Dynkin diagram
of $D_4^{++}$, the hyperbolic over-extension of the $D_4$ Dynkin
diagram, is depicted in  fig.~\ref{d4pp}. Recall that the 
Weyl group of $D_4$ is $W(D_4) = 2^3\rtimes \Sg_4$. 

\begin{figure}
\centering
\begin{picture}(110,80)
\thicklines
\multiput(10,40)(30,0){4}{\circle*{8}}
\put(70,70){\circle*{8}}
\put(70,10){\circle*{8}}
\put(14,40){\line(1,0){82}}
\put(70,66){\line(0,-1){52}}
\put(4,25){$-1$}
\put(37,25){$0$}
\put(60,25){$1$}
\put(97,25){$3$}
\put(82,68){$2$}
\put(82,8){$4$}
\end{picture}
\caption{\label{d4pp}\sl Dynkin diagram of $D_4^{++}$ with numbering
 of nodes.}
\end{figure}
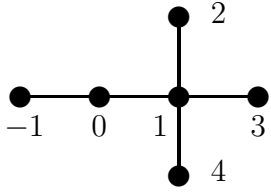

Our central result is 
\begin{prop} $W^+(D_4^{++})\cong \PGL$.
\end{prop}

\noindent{\bf{Proof:}}
Substituting (\ref{D4roots}) into the expressions of the generating matrices 
$S_i$, $0\leq i\leq 4$, given in (\ref{W+}), one sees that all these matrices 
have ${\text{det}}(S_i S_i^\dagger) = 1$ and $\rmD(S_i) = 1$ 
so they belong to  $\PGL\subset PSL_2(\Hq)$ as defined in
(\ref{gq}). Therefore, we get 
$W^+(D_4^{++}) \subset \PGL$. To prove the converse we show that all the
generating matrices of Proposition~\ref{psl0gen} are in
$W^+(D_4^{++})$. To do this  
we compute  
\be
S_i S_1^{-1} &=& \lp\begin{array}{cc}\ve_i&0\\0&\bve_i\end{array}\rp
\qquad \mbox{for $i=2,3,4$} \; ,\nn\\
S_2 S_3 &=& \lp\begin{array}{cc}i&0\\0&j \end{array}\rp\; , \;\;
S_3 S_4 = \lp\begin{array}{cc}j&0\\0&k \end{array}\rp\;,\;\;
S_4 S_2 = \lp\begin{array}{cc}k&0\\0&i \end{array}\rp\,.
\ee
It is straightforward to see that by further multiplication and 
permutation we can obtain from these matrices and $S_1$ any matrix of the form
\be\label{finiteD4}
\lp\begin{array}{cc} a&0\\0&b\end{array}\rp \quad \mbox{and} \quad
\lp\begin{array}{cc} 0&a\\b&0\end{array}\rp 
\ee
where $(a,b)$ is any pair of Hurwitz units whose product is a Lipschitz unit, 
i.e. $N(a) = N(b) = 1$ and $ab\equiv 1\ (\text{mod}\ \Cq)$. We note that up to
$(a,b)\cong(-a,-b)$ there are $96$ pairs with this property so that there are
$96$ matrices of each type in (\ref{finiteD4}) in $PSL_2(\Hq)$ corresponding
to the even and odd part of the finite Weyl group $W(D_4)$ of order $192$.  

Given the matrices (\ref{finiteD4}) we can
also reconstruct the translation matrix from $S_0$ (using $\theta^3=-1$) and
\be\label{transd4}
\lp\begin{array}{cc} 0& \theta^2\\ \bth^2 &0\end{array}\rp \cdot 
\lp\begin{array}{cc} 0& \theta\\ - \bth &1\end{array}\rp  \cdot
\lp\begin{array}{cc} - \bth & 0\\ 0 & \theta\end{array}\rp = 
\lp\begin{array}{cc} 1&-1\\ 0&1\end{array}\rp\,.
\ee 
The translation matrix is the inverse of this, giving the last of the 
required generators (\ref{gl2hgen}). $\square$\medskip

The appearance of the group $\PGL$ for $D_4^{++}$ (rather than $PSL_2(\Hq)$) 
may be viewed as a manifestation of the triality symmetry of the $D_4$
Dynkin diagram. Since the even cyclic permutation of $(i,j,k)$ is
realized by conjugation by the Hurwitz unit $\theta$, the associated
diagonal matrix $\text{diag}(\theta,\theta)\in PSL_2(\Hq)$ solves the
necessary conditions (\ref{rtsym}), but itself is not part of the Weyl 
group, and therefore has to be removed from $PSL_2(\Hq)$. This is 
what the $\rmD S =1$ 
condition achieves, and since $\theta^3=-1$
we find that the subgroup $\PGL$ is index three in $PSL_2(\Hq)$, in
accordance with Lemma~\ref{d4lemma}. In this way, the group of (even)
{\em outer} automorphism of $D_4$ can also be realized by matrix
conjugation. 

For later reference we denote the alternative generating set of 
the group $\PGL$ furnished by the even $D_4^{++}$ Weyl group by
\be\label{Ggen}
g_i = S_i^{D_4} \quad\quad \text{for}\,\,i=0,\ldots,4\,.
\ee

An alternative argument leading to the statement of the proposition 
can be based on the fact that {\em every} Hermitian matrix
$\capX\in\Lambda(D_4)$  
obeying $\det \capX= -1$ is a real root,\footnote{It is here that the 
hyperbolicity of the Kac--Moody algebra is essential; this statement is 
no longer true if the indefinite Kac--Moody algebra is not hyperbolic.}
and that the (full) Weyl group acts transitively on the set of real roots.
In the case of $D_4^{++}$ they form  a single Weyl orbit, since the Dynkin
diagram has only single lines. This implies that we can generate {\em all}
real roots from the hyperbolic simple root $\alpha_{-1}$ by acting with the
full Weyl group, {\em viz.} 
\be
\Delta^{\text{re}} = \big\{ \capX\in \Lambda(D_4) \,| \,\det\capX = -1\big\}
           = W\cdot\lp\begin{array}{cc}1&0\\0&-1\end{array}\rp\,.
\ee
In particular, considering all even Weyl images of the form 
$s(\alpha_{-1}) = S\alpha_{-1} S^\dagger$, and using the fact that
Weyl transformations preserve the norm, we get
\be
 -1= \det (\alpha_{-1}) = \det(S\alpha_{-1}S^\dagger) = 
  - \det(SS^\dagger)
\ee
where the third equality can be verified by direct computation for
the quaternions. This is the same condition as (\ref{PSLH}). Conversely,
to any $S\in PSL_2(\Hq)$ satisfying  (\ref{PSLH}) we can associate the real root
$S\alpha_{-1}S^\dagger$, whence $S$ becomes associated to some symmetry of
$\Lambda(D_4)$. This seems to give all of $PSL_2(\Hq)$ but the 
reasoning does not distinguish between inner and outer transformations. As
argued above, taking the outer (diagram) automorphisms into account is the
same as demanding the extra condition (\ref{gq}).

\end{subsection}

\begin{subsection}{Even Weyl group $W^+(B_4^{++})$}

The $B_4$ root lattice is isomorphic to the hypercubic lattice of the 
Lipschitz integers (in the same way that $B_2$ was associated with the 
cubic lattice of Gaussian integers). The simple roots can be chosen 
as follows
\begin{align}\label{typeB4}
\ve_1 &= a_1 = 1 \;\; , &
\ve_2 &= a_2 = \frac12 \big( -1 + i - j - k \big) \; , \nn\\
\ve_3 &= a_3 =\frac12 \big( -1 - i + j - k \big) \; , &
\ve_4 &= \sqrt{2} a_4 = \frac{-j+k}{\sqrt{2}} \; .
\end{align}
The highest (long) root is
\be
\theta = a_2 + 2 a_1 + 2 a_3 + 2 a_4 = \frac12 \big( 1 - i - j - k \big)\,.
\ee
Hence, $a_4$ is the short simple root, while the other (long) roots,
including the highest root $\theta$, are normalized to unity. These
simple roots were chosen to be close to the $D_4$ simple roots: $a_1$,
$a_2$ and $a_3$ 
agree and only $a_4$ is changed; the highest roots agree. This 
somewhat obscures the fact that the $B_4$ lattice is a scaled and rotated 
version of the lattice of Lipschitz numbers but makes the determination 
of the Weyl group simpler. The finite 
Weyl group is $W(B_4) = 2^4\rtimes \Sg_4$, and is concretely realized 
via permutations and reflections of the four basis elements $1,i,j,k$. 
Observe also that the above set of simple roots is not invariant under
cyclic (triality) transformations of the three imaginary units, as was
to be expected.

\begin{figure}
\centering
\begin{picture}(140,80)
\thicklines
\multiput(10,40)(30,0){5}{\circle*{8}}
\put(10,40){\line(1,0){90}}
\put(70,70){\circle*{8}}
\put(70,66){\line(0,-1){26}}
\put(100,44){\line(1,0){30}}
\put(100,36){\line(1,0){30}}
\put(110,50){\line(1,-1){10}}
\put(110,30){\line(1,1){10}}
\put(4,25){$-1$}
\put(37,25){$0$}
\put(67,25){$1$}
\put(97,25){$3$}
\put(127,25){$4$}
\put(82,68){$2$}
\end{picture}
\caption{\label{b4pp}\sl Dynkin diagram of $B_4^{++}$ with numbering
 of nodes.}
\end{figure}
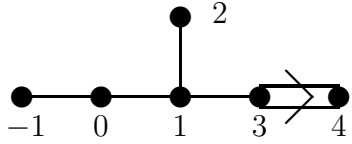

The hyperbolic extension $B_4^{++}$ can be obtained by following the
procedure described before and leads to the Dynkin diagram of fig.~\ref{b4pp}. 
Turning to its even hyperbolic Weyl group $W^+(B_4^{++})$ we notice that 
among the generators $S_i\in PSL_2(\quat)$ (for $i=0,\ldots,4$) the 
first four are identical to those of  $W^+(D_4^{++})$. By explicit 
computation one finds furthermore that
\be\label{b4g4}
S_4S_3S_4  
= \lp\begin{array}{cc}
  0&\frac12(1+i+j-k)\\
   \frac12(- 1+i+j-k)&0\end{array}\rp =  g_4 \,,
\ee
where $g_4$ was defined in (\ref{Ggen}) in the alternative generating
set of $\PGL$.
This shows that $W^+(B_4^{++})$ has a subgroup isomorphic to $\PGL$
generated by $S_0$, $S_1$, $S_2$, $S_3$ and $S_4S_3S_4$. Adjoining to 
these the element $S_4$ which satisfies $S_4^2= - \id \cong \id$ 
gives an extension in which $\PGL$ is a normal subgroup of index 2, 
and we arrive at

\begin{prop}
$ W^+(B_4^{++})\cong \PGL \, \rtimes\, 2$.
\end{prop}

The semi-directness follows from expanding the action $S_4 S S_4^{-1}$
for $S\in \PGL$ and using the properties of $\ve_4$. As the index $2$
extension is realized via the matrix $S_4$ which contains $1/\sqrt{2}$ 
in all entries we see that $W^+(B_4^{++})$ consists of all matrices 
(\ref{psl2hmat}) in $PSL_2(\quat)$ such that either all $a,b,c,d$ are 
Hurwitz numbers, or all are pure octahedral numbers $\in i_{\text{O}}\cdot \Hq$ 
(it is easy to see that this property is preserved under matrix 
multiplication because $i_O\Hq \cdot \Hq = i_{\text{O}}\Hq$ and
$i_{\text{O}}\Hq \cdot i_{\text{O}}\Hq = \Hq$). This is the `trick' by which the 
discreteness of the matrix group is achieved despite the 
fact that $\cR$ is dense in $\quat$.

To conclude the analysis of $B_4^{++}$ we relate our results, 
restricted to the finite Weyl group $W^+(B_4)$, to the classification of
\cite{CoSm03} (see their Table~4.2) where $W^+(B_4)$ is denoted by $\pm\frac16
[O\times O]$. It is not hard to check that all the generating elements given
in \cite{CoSm03} can be obtained from products of matrices of the type
(\ref{SiSj}) when the $B_4$ units (\ref{typeB4}) are plugged in.

\end{subsection}

\begin{subsection}{Even Weyl group $W^+(C_4^{++})$}

For $C_4$ the direction of the arrow in the finite Dynkin diagram is 
reversed compared to $B_4$. A convenient choice of simple roots requires again
quaternions which are not Hurwitz. In order to facilitate the comparison
with $D_4$ we choose the first three simple roots of $C_4$ to coincide
with the simple roots $a_3, a_1, a_4$ of $D_4$ up to a re-scaling
required to maintain the unit normalization of the highest root 
$\theta$. Thus
\begin{align}\label{typeC4}
\ve_1 &= \sqrt{2} a_1 = 1 \;\; , & 
\ve_2 &= \sqrt{2} a_2 = \frac12 \big( -1 + i - j - k \big) \nn\\
\ve_3 &= \sqrt{2} a_3 =\frac12 \big( -1 - i + j - k \big) \;\; , &
\ve_4 &= a_4 = \frac{-j+k}{\sqrt{2}}
\end{align}
with the (long) highest root
\be\label{C4hst}
\theta = 2 a_1 + 2 a_2 + 2 a_3 + a_4 = \frac{ - j - k }{\sqrt{2}} \,,
\ee

\begin{figure}
\centering
\begin{picture}(170,80)
\thicklines
\multiput(10,40)(30,0){6}{\circle*{8}}
\put(70,40){\line(1,0){60}}
\put(130,44){\line(1,0){30}}
\put(130,36){\line(1,0){30}}
\put(150,50){\line(-1,-1){10}}
\put(150,30){\line(-1,1){10}}
\put(10,40){\line(1,0){30}}
\put(40,44){\line(1,0){30}}
\put(40,36){\line(1,0){30}}
\put(50,50){\line(1,-1){10}}
\put(50,30){\line(1,1){10}}
\put(4,25){$-1$}
\put(37,25){$0$}
\put(67,25){$2$}
\put(97,25){$1$}
\put(127,25){$3$}
\put(157,25){$4$}
\end{picture}
\caption{\label{c4pp}\sl Dynkin diagram of $C_4^{++}$ with numbering
 of nodes.}
\end{figure}
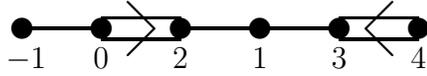

For the determination of the Weyl group we again choose a convenient basis.
The hyperbolic extension $C_4^{++}$ has Dynkin diagram given in 
fig.~\ref{c4pp}. Although the $B_4$ and $C_4$ Weyl groups are isomorphic 
(cf. table~\ref{integers}), the Weyl groups of $C_4^{++}$ and $B_4^{++}$ 
are nevertheless different. This difference in our basis is reflected
only in the difference between the highest roots; for $C_4$ we have an
octahedral unit of order four whereas for $B_4$ the highest root was a
Hurwitz number of order six. 

Proceeding as for $W^+(B_4^{++})$ we first observe that all $S_i$ for 
$i=1,\ldots,4$ are identical to those of $B_4^{++}$ (since the finite 
Weyl groups are isomorphic), and hence again 
\be\label{c4g4}
S_4S_3S_4  
= \lp\begin{array}{cc}
  0&\frac12(1+i+j-k)\\ \frac12(-1+i+j-k)&0\end{array}\rp =  g_4 \,
\ee
in terms of (\ref{Ggen}). Therefore we obtain all diagonal and
off-diagonal matrices 
(\ref{finiteD4}) reflecting the fact that the finite Weyl groups of
$B_4$ and $C_4$ are isomorphic. Using the $C_4$ highest root
(\ref{C4hst}) we now define 
\be\label{twgs}
\tilde{g}_i = \lp\begin{array}{cc}   \theta&0\\     0& 1\end{array}\rp g_i
\lp\begin{array}{cc}   \bar\theta&0\\     0& 1\end{array}\rp\quad\text{for}
\;\;i=1,2,3,4\,,
\ee
which is a unitary transformation of all the generators of $\PGL$ defined 
in (\ref{Ggen}) except for $g_0$. The matrices $\tilde{g}_i$ ($i=1,2,3,4$) 
also belong to $W^+(C_4^{++})$ since $\theta a \bar\theta b \equiv ab\,
(\text{mod}\,\Cq)$. Furthermore, since
\be
\lp\begin{array}{cc}   \theta& 0\\     0& 1\end{array}\rp
\lp\begin{array}{cc}   1&1\\     0& 1\end{array}\rp
\lp\begin{array}{cc}   \bar\theta& 0\\     0& 1\end{array}\rp
=\lp\begin{array}{cc}   1& \theta\\     0& 1\end{array}\rp\,.
\ee
we also have the unitarily transformed translation matrix. From this
we conclude that the group generated by $S_0$, $S_1$, $S_2$, and $S_3$
and $S_4S_3S_4$ consists of matrices of the type 
\be
\lp\begin{array}{cc}   \theta& 0\\     0& 1\end{array}\rp 
\lp\begin{array}{cc}   a& b\\     c& d\end{array}\rp
\lp\begin{array}{cc}   \bar\theta& 0\\     0& 1\end{array}\rp\,\quad
\text{for}\, \lp\begin{array}{cc}   a& b\\     c& d\end{array}\rp\in \PGL\,.
\ee
We denote this transformed (or `twisted') $\PGL$ by $\tPGL$. In order to 
describe the full $W^+(C_4^{++})$ we still need to adjoin the matrix 
$S_4$ which still satisfies $S_4^2=-\id\cong \id$ so that we find

\begin{prop}
$W^+(C_4^{++})  \cong \tPGL\rtimes 2\,.$
\end{prop}

The semi-directness of the product can be verified by expanding the action of
$S_4$ on $\tPGL$ which can be seen to give an automorphism of $\tPGL$.  
Although $\tPGL\cong \PGL$, it must be stressed that the Weyl groups of 
$C_4^{++}$ and $B_4^{++}$ are nevertheless different, since for $C_4^{++}$ 
the extension $S_4$ matrix is not unitarily transformed, so that the 
group extension is different.\footnote{This is similar to the affine 
  Weyl groups $W^+(B_4^+)$ and $W^+(C_4^+)$ although both are semi-direct 
  products of isomorphic finite groups with abelian groups of the same rank. 
  However, the automorphism involved in the definition of the semi-direct 
  product is different.}
The $PSL_2(\quat)$ matrices (\ref{psl2hmat}) belonging to
$W^+(C_4^{++})$ can be described by saying that either  
$a,d\in\Hq$ and $b,c\in i_{\text{O}}\cdot\Hq$, or conversely 
$a,d\in i_{\text{O}}\cdot\Hq$ and $b,c\in\Hq$. Again, one can easily verify 
that this structure is preserved under multiplication.

Since the finite Weyl group satisfies $W(C_4)\cong W(B_4)$, the relation to the
chiral group $\pm\frac16 [O\times O]$ of the classification in \cite{CoSm03}
holds analogously to the last section.

\end{subsection}

\begin{subsection}{Even Weyl group $W^+(F_4^{++})$}

The root system of $F_4$ consists of two copies of the $D_4$ system, 
one of which is rescaled. There are two short simple roots and two
long simple roots. We choose the simple roots as follows 
\begin{align}\label{typeF4}
\ve_1 &= a_1 = 1 \;\; , & 
\ve_2 &= a_2 = \frac12 \big( -1 + i - j - k \big) \; ,\nn\\
\ve_3 &= \sqrt{2} a_3 = \frac{-i+j}{\sqrt{2}}\;\; , &
\ve_4 &= \sqrt{2} a_4 =\frac{ -j+k}{\sqrt{2}} \; ;
\end{align}
the highest (long) root is 
\be
\theta = 2 a_1 + 3 a_2 + 4 a_3 + 2 a_4 = \frac12 \big( 1 - i - j - k \big)\,.
\ee
The hyperbolic extension  $F_4^{++}$ has the Dynkin diagram 
shown in fig.~\ref{f4pp}.

\begin{figure}
\centering
\begin{picture}(170,80)
\thicklines
\multiput(10,40)(30,0){6}{\circle*{8}}
\put(10,40){\line(1,0){90}}
\put(100,44){\line(1,0){30}}
\put(100,36){\line(1,0){30}}
\put(110,50){\line(1,-1){10}}
\put(110,30){\line(1,1){10}}
\put(130,40){\line(1,0){30}}
\put(4,25){$-1$}
\put(37,25){$0$}
\put(67,25){$1$}
\put(97,25){$2$}
\put(127,25){$3$}
\put(157,25){$4$}
\end{picture}
\caption{\label{f4pp}\sl Dynkin diagram of $F_4^{++}$ with numbering
 of nodes.}
\end{figure}
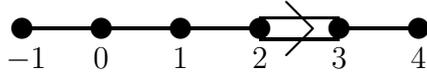

The even Weyl group $W^+(F_4^{++})$ can be related to $\PGL$ in the 
following way: The generators $g_0$, $g_1$ and $g_2$ are identical 
to $S_0$, $S_1$ and $S_2$ constructed from the $F_4$ simple roots  
(\ref{typeF4}). Now consider
\be\label{f4g3}
S_3S_2S_3
= \lp\begin{array}{cc}
  0&\frac12(-1-i+j-k)\\
   \frac12(1-i+j-k)&0\end{array}\rp
\equiv g_3 
\ee
and
\be\label{f4g4}
S_4S_3S_2S_3S_4
= \lp\begin{array}{cc}
  0&\frac12(-1-i-j+k)\\
   \frac12(1-i-j+k)&0\end{array}\rp
\equiv g_4\,, 
\ee
in terms of (\ref{Ggen}),
so that we find all $\PGL$ generators in $W^+(F_4^{++})$. Adjoining 
to $\PGL$ the generators $S_3$ and $S_4$ one obtains all of $W^+(F_4^{++})$. 
Since $S_3$ and $S_4$ together generate the symmetric group on three 
letters we arrive at
\begin{prop}
$W^+(F_4^{++}) \cong \PGL \, \rtimes \, \Sg_3\cong PSL_2(\Hq)\,\rtimes\, 2$.
\end{prop}
The last equality follows because the (diagonal) entries of $S_3S_4$
are pure Hurwitz units, and both of order three, generating the cyclic
group  $\ints_3$ which extends $\PGL$ to $PSL_2(\Hq)$. The
semi-directness of the product can be verified explicitly by using the
properties of $\ve_3$ and $\ve_4$. 

The finite Weyl group $W^+(F_4)$ is now the chiral group $\pm\frac12 [O\times
O]$ in the classification of \cite{CoSm03}. Again it is not hard to see that
the units (\ref{typeF4}) give rise to the generating elements given there.

\end{subsection}

\begin{subsection}{Even Weyl group $W^+(A_4^{++})$}

We discuss the algebra $A_4^{++}$ last, because, somewhat surprisingly,
it is the most cumbersome to describe in the present framework among 
the rank 6 algebras. In order
to embed its root system into the quaternions one needs to make use of the 
{\em icosians} $\Iq$ introduced at the end of Section~5.2. The algebra
$A_4^{++}$ is simply laced, and we choose the following four units
\begin{align}\label{typeA4}
\ve_1 &= a_1 = 1\,,&\quad
\ve_2 &= a_2 = \frac12 ( -1 +i -j - k    )\,,&\nn\\
\ve_3 &= a_3 = -i \,,&\quad
\ve_4 &= a_4 = \frac12 (i + \si j  + \tau k )\,,&
\end{align}
which have inner products corresponding to the $A_4$ Cartan matrix.
The highest root is
\be
\theta = a_1 + a_2 + a_3 + a_4 = \frac12 (1 - \tau j - \si k)\,.
\ee
The Dynkin diagram of $A_4^{++}$ is displayed in fig.~\ref{a4pp}. 

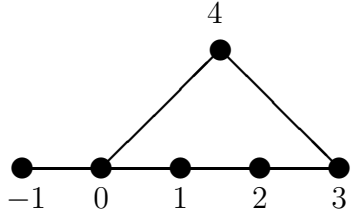
\begin{figure}
\centering
\begin{picture}(170,100)
\thicklines
\multiput(10,40)(30,0){5}{\circle*{8}}
\put(85,85){\circle*{8}}
\put(10,40){\line(1,0){120}}
\put(100,40){\line(1,0){30}}
\put(100,40){\line(1,0){30}}
\put(40,40){\line(1,1){45}}
\put(85,85){\line(1,-1){45}}
\put(4,25){$-1$}
\put(37,25){$0$}
\put(67,25){$1$}
\put(97,25){$2$}
\put(127,25){$3$}
\put(80,95){$4$}
\end{picture}
\caption{\label{a4pp}\sl Dynkin diagram of $A_4^{++}$ with numbering
 of nodes.}
\end{figure}

As a first step we identify the finite Weyl group $W(A_4)$ of order 
$120$ in terms of diagonal and off-diagonal matrices. The diagonal matrices 
can be generated from $\text{diag}\,(\ve_i,\bar{\ve}_i)$ for $i=2,3,4$ and, 
modulo $-\id$, one obtains all $60$ even elements of $W(A_4)$ in 
this way. This shows that the root system of $A_4$ is not multiplicatively 
closed. Closer inspection reveals that the diagonal and off-diagonal 
matrices are of the form
\be\label{finiteA4}
\lp\begin{array}{cc}   a& 0\\     0& a^*\end{array}\rp\quad\text{and}\quad
\lp\begin{array}{cc}   0& a\\     -a^*& 0\end{array}\rp\,,
\ee
where
\be\label{autoa4}
a^* = a'_0  -  a'_1 i - a'_3 j - a'_2 k 
  \quad\quad \text{for}\,\, a=a_0 +a_1 i +a_2 j +a_3 k
\ee
with $(p + \sqrt{5} q)' := p-\sqrt{5} q$ for $p,q\in \rats$. 
For icosians the conjugation in $\sqrt{5}$ amounts to interchanging 
$\tau$ and $\sigma$ everywhere. The operation (\ref{autoa4}) is an
involutive automorphism of the ring $\cI$ of icosians and satisfies
in  particular $(ab)^*=a^* b^*$. The automorphism is outer due to the
`conjugation' in $\sqrt{5}$; however, the exchange of $j$ and $k$ together 
with the complex conjugation alone is an {\em inner automorphism} of the
quaternion algebra $\quat$, and is explicitly realized by conjugation 
with the octahedral unit $i_{\text{O}}$ from (\ref{ounits}).   

Since there are $120$ icosian units there 
are $60$ matrices of either type in (\ref{finiteA4}) after quotienting 
out by $\{\id,-\id\}$ in agreement with the order of the Weyl group
$W(A_4)\cong \Sg_5$. Moreover, 
the $A_4$ root lattice is defined by the condition  
\be\label{A4latt}
Q_{A_4} = \left\{ z \in \cI\;\vline\; z^*=\bar{z}\right\}\,.
\ee
Because conjugation reverses the order of factors but the
$*$-automorphism does not, we see from (\ref{A4latt}) again that
$Q_{A_4}$ is not multiplicatively closed and therefore does not form a
subring of $\cI$. $Q_{A_4}$ is closed under complex
conjugation.\footnote{The fixed point set of $*$ on the set of icosian
 units is $\{1\}$.} The action of the even finite Weyl group on $z\in
Q_{A_4}$ is computed to be $z\mapsto a z\bar{a}^*$, which manifestly
preserves the root lattice. 

The full $W^+(A_4^{++})$ is obtained by adjoining to (\ref{finiteA4})
the matrix 
\be\label{a4t}
\lp\begin{array}{cc}   1& \theta\\     0& 1\end{array}\rp
\ee
and the result is
\begin{prop}
$W^+(A_4^{++}) \cong PSL_2^{(0)}(\cI)$.
\end{prop}
Here, $PSL_2^{(0)}(\cI)$ is defined as the discrete subgroup of
$PSL_2(\Iq)$  satisfying the conditions (\ref{rtsym}) with the four
units $a_i$ from (\ref{typeA4}); the
$(0)$-superscript indicates that additionally the $\ints_2$ outer
automorphism of the $A_4$ Dynkin diagram has to be removed. 
This outer automorphism acts concretely as $z\mapsto -z$ on $z\in
Q_{A_4}$ and the corresponding matrix is $\text{diag}(1,-1)$ which has
to be quotiented out. The explicit result of solving the conditions
(\ref{rtsym}) and taking the quotient for diagonal and
off-diagonal matrices leaves the ones listed in (\ref{finiteA4}). For
the general matrices $S$ it leaves in particular the matrix (\ref{a4t}).

The finite Weyl group $W^+(A_4)$ in terms of the classification of
\cite{CoSm03} is $+\frac1{60}[I\times\bar{I}]$ and again we see that there the
generators given there can be reproduced from the icosian units (\ref{typeA4}).

One might hope to give a more natural description of $PSL_2^{(0)}(\cI)$ as the
kernel of a group homomorphism analogous to $\rmD$ used to define
$PSL_2^{(0)}(\Hq)$, but since the ideal $\Cq_\cI$ generated by the
commutators in $\cI$ is not proper, this idea fails. 

\end{subsection}

\begin{subsection}{Icosians and $E_8$}

Using icosians it is also possible to give an embedding of the $E_8$ root
system into $\quat$ at the cost of introducing a new inner product. This
inner product is defined by 
\be\label{icoE8}
\left( a, b \right)_\tau := x\,,\quad \text{if}\,\, (a,b)= a\bar{b}+b\bar{a} =
x+\tau y
\ee
for $a,b,x,y$ being rational quaternions. The $E_8$ simple roots are
then given by (see 
e.g. \cite{MoPa93}) 
\begin{align}\label{e8quat}
a_1 &= \frac12 ( -\si -\tau i            -k      )\,,&\quad
a_2 &= \frac12 (     -\si i   -\tau j    +k      )\,,&\nn\\
a_3 &= \frac12 (     +i      -\si j    -\tau k   )\,,&\quad
a_4 &= \frac12 (     -\tau i   +j       +\tau^2 k )\,,&\nn\\
a_5 &= \frac12 (     +\tau i   +j       -\tau^2 k )\,,&\quad
a_6 &= \frac12 (     +i      -\tau^2 j  +\tau k   )\,,&\nn\\
a_7 &= \frac12 ( +1  -\tau^2 i          -\tau k   )\,,&\quad
a_8 &= \frac12 (     -i       -\si j   +\tau k   )\,.&
\end{align}
Note that $\tau^2=\tau+1$. Here, $a_4$, $a_5$, $a_6$ and $a_7$ belong to
$\tau\Iq$ constituting the second half of the $240$ $E_8$ roots. The 
subspace spanned by integer linear combinations of the basis
vectors (\ref{e8quat}) is dense in $\quat$ (as one would expect, since 
we are projecting an 8-dimensional lattice onto a 4-dimensional hyperplane 
in such a way that distincts points remain distinct),
and the same is, of course, true for the ring generated by the above 
elements in $\quat$. In the following section we will present and discuss 
an {\em octonionic} realization of the $E_8$ root system, and we will see 
that the non-associativity leads to drastic changes in the matrix 
realization of the Weyl group. For this reason, one might have hoped 
to be able to avoid problems with non-associativity by working with the
above quaternionic realization, but unfortunately the difficulties remain,
mainly due to the modified inner product (\ref{icoE8}). Namely, the 
projection orthogonal to $\tau$, which enters (\ref{icoE8}), invalidates 
our Theorem~\ref{Weylprop}, and we therefore cannot use the above 
representation to find a matrix representation of the $E_8$ Weyl group.
For this reason we do not pursue this quaternionic realization of $E_8$ 
further in this paper.

\end{subsection}

\end{section}

\begin{section}{Octonions $\Kk=\oct$}
\label{octsec}

In this section we turn to the non-associative octonions and their relation to
hyperbolic algebras of maximal rank $10$. We restrict our attention mostly to
the case $E_8^{++}\equiv E_{10}$ and exhibit in this case the kind of new
complications arising due to the non-associativity of $\oct$. For the other
two hyperbolic over-extended algebras $B_8^{++}$ and $D_8^{++}$ we content
ourselves with briefly presenting the octonionic realizations of their root
lattices and stating their relation to that of $E_{10}$.

\begin{subsection}{Octavians and the $E_8$ lattice}

We use the octonionic multiplication conventions of Coxeter
\cite{Co46}. Hence, the seven imaginary units $e_i$ appearing in the
expansion of an octonion $z$ via ($e_0=1$)
\be
z= n_0 +\sum_{i=1}^7 n_i e_i
\ee
multiply (in quaternionic subalgebras) according to
\be
e_i e_{i+1} e_{i+3} = -1\,,
\ee
where the indices are to be taken modulo seven and the relation is
totally antisymmetric.
The multiplication table for the octonions can be neatly summarized by
using the Fano plane, as shown in Fig.~\ref{fanofig}. The
following notation for the structure constants of the imaginary units
can also be used 
\be
e_i e_j =-\delta_{ij}+ f_{ijk} e_k
\ee
where we adopt the usual summation convention from now on. The
structure constants $f_{ijk}$ are totally antisymmetric, and satisfy
(see e.g. \cite{deWiNi84} where many further identities involving these
structure constants can be found)
\be
f_{ijk}f^{kmn} = 2 \delta^{mn}_{ij} -\frac16\epsilon_{ij}{}^{mnrst}f_{rst}\,,
\ee
with $\delta^{mn}_{ij} := \frac12\delta^m_i\delta^n_j -
\frac12\delta^n_i\delta^m_j$.

\begin{figure}[t!]
\centering\includegraphics[scale=.9]{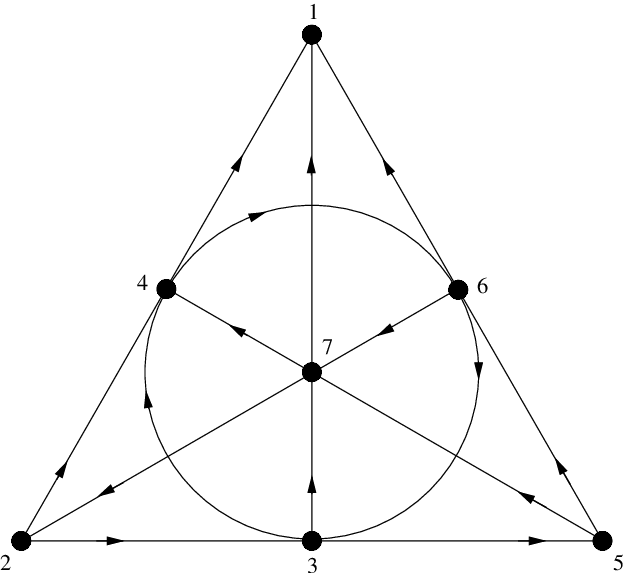}
\caption{\label{fanofig} \sl The Fano plane illustrating associative
 triples and octonionic multiplication.}
\end{figure}

We introduce a few concepts which are useful for studying the $E_8$ and
$E_8^{++}$ Weyl groups. 
The {\em associator} of octonions is defined by
\be
\left[ a, b, c \right] = a (bc) -(ab) c
\ee
and is totally antisymmetric in $a,b,c$ by virtue of the 
alternativity of $\oct$. Put differently, any subalgebra generated by two
elements is associative. In components it reads
\be
[a,b,c]= [a,b,c]_i e_i \,\quad\text{with}\quad 
[a,b,c]_i = -\frac13 \epsilon_i{}^{jklrst}a_jb_kc_l f_{rst}\,.
\ee

The Jacobi identity gets modified by an associator term to
\be
\left[\left[ a, b\right], c\right] 
 +\left[\left[ b, c\right], a\right]
 +\left[\left[ c, a\right], b\right]
 +6 \left[a, b, c\right] = 0\,.
\ee 

The {\em Moufang laws} can be written as
\be\label{mou}
x (yz) = (xyx) (x^{-1}z)\,,\quad\quad
(yz) x = (yx^{-1}) (xzx)\,.
\ee 
They imply for example that $a(xy)a=(ax)(ya)$. Similarly one can transform
$axa^{-1}$ to show that this is an automorphism whenever $a\in\oct$ satisfies $a^3=\pm 1$.
Another consequence of alternativity is that the following expressions are well defined without brackets
\be
a x a = (ax) a= a (xa)\,,\quad\quad
a x a^{-1} = (ax) a^{-1} = a (xa^{-1})\,.
\ee
The first one implies the second one since $a^{-1} = \bar{a}/N(a)$. These
relations are the basic reason why our formulas (\ref{MWeyl}) with the
generating matrices (\ref{Weyl2}) also make sense for octonions.

Coxeter (following Bruck) has established that there is a maximal order of
integers within $\oct$ \cite{Co46}. We briefly recall the construction from
\cite{CoSm03}. One starts from the sets of all Hurwitz numbers for all
associative triples $(1,e_i,e_{i+1},e_{i+3})$ of fig.~\ref{fanofig} and their
complements, arriving at the `Kirmse numbers'. This is not a maximal order 
but will be so after exchanging $1$ with any of the seven imaginary units; 
for definiteness we choose this unit as $e_3$. In addition, the lattice 
of integers is a scaled copy of the $E_8$ root lattice. The following 
octonionic units correspond to the simple roots of $E_8$ (labeled as 
in fig.~\ref{e8pp}) 
\be\label{E8simple}
a_1 = e_3 \quad &,& \quad 
a_2 = \frac12\big( -e_1 -e_2 - e_3 + e_4 \big) \nn\\
a_3 = e_1 \quad &,& \quad 
a_4 = \frac12\big( - 1 -e_1 - e_4 + e_5 \big) \nn\\
a_5 = 1 \quad &,& \quad 
a_6 = \frac12\big( - 1 -e_5 - e_6 - e_7 \big) \nn\\
a_7 = e_6 \quad &,& \quad 
a_8 = \frac12\big( - 1 + e_2 + e_4 + e_7 \big) \,.
\ee
Their integer span gives all octonionic integers and, following \cite{CoSm03},
we call these integers {\em octavians} and denote them by $\Oo$.
The highest $E_8$ root is represented by the unit
\be
\theta &=& 2a_1 + 3 a_2 + 4 a_3 + 5 a_4 + 6 a_5 + 4 a_6 + 2 a_7 + 3 a_8 \nn\\
      &=& \frac12 \big( e_3 + e_4 + e_5 - e_7 \big) \,.
\ee
There is a total of $240$ unit octavians corresponding to the $240$ roots of
$E_8$.  These definitions imply that no octavian contains only a quaternionic triple or its complement in its expansion.

\begin{figure}
\centering
\begin{picture}(260,80)
\thicklines
\multiput(10,40)(30,0){9}{\circle*{8}}
\put(190,70){\circle*{8}}
\put(14,40){\line(1,0){232}}
\put(190,44){\line(0,1){22}}
\put(4,25){$-1$}
\put(37,25){$0$}
\put(67,25){$1$}
\put(97,25){$2$}
\put(127,25){$3$}
\put(157,25){$4$}
\put(187,25){$5$}
\put(217,25){$6$}
\put(247,25){$7$}
\put(202,68){$8$}
\end{picture}
\caption{\label{e8pp}\sl Dynkin diagram of $E_{10}\equiv E_8^{++}$
 with numbering of nodes.}
\end{figure}
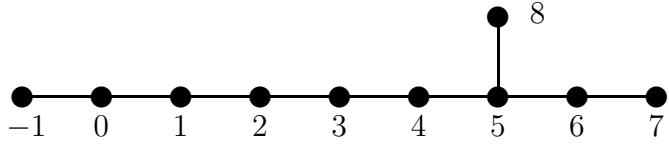

\end{subsection}

\begin{subsection}{$D_8$ and $B_8$ root systems}

Before studying the $E_8$ and $E_{10}$ Weyl groups we give the octonionic
realizations of the root lattices of $D_8$ and $B_8$ 
which can be extended to the hyperbolic over-extensions in the usual manner.

For $D_8$ the simple roots are (as labeled in fig.~\ref{d8pp})
\be\label{D8simple}
D_8:\quad\quad\quad a_1 = e_3 \quad &,& \quad 
a_2 = \frac12\big( -e_1 -e_2 - e_3 + e_4 \big) \nn\\
a_3 = e_1 \quad &,& \quad 
a_4 = \frac12\big( - 1 -e_1 - e_4 + e_5 \big) \nn\\
a_5 = 1 \quad &,& \quad 
a_6 = \frac12\big( - 1 - e_5 - e_6 - e_7 \big) \nn\\
a_7 = \frac12\big( e_2 - e_3 +e_6 -e_7 \big)\quad &,& \quad 
a_8 = \frac12\big( - 1 + e_2 + e_4 + e_7 \big) \,.
\ee
Compared to the $E_8$ simple roots of (\ref{E8simple}) the only difference is
in the simple root $a_7$ which is a specific linear combination describing an
embedding of $D_8$ into $E_8$ \cite{SchnWe04,KlNi04}. The highest root of
$D_8$ is given by 
\be
\theta^D &=& 2a_1 + 2 a_2 + 2 a_3 + 2 a_4 + 2 a_5 +  a_6 +  a_7 +  a_8 \nn\\
      &=& \frac12 \big( e_3 + e_4 + e_5 - e_7 \big) \,.
\ee
which is identical to the highest $E_8$ root. This will make it easy to relate
the two hyperbolic Weyl groups.

\begin{figure}
\centering
\begin{picture}(260,80)
\thicklines
\multiput(10,40)(30,0){8}{\circle*{8}}
\put(190,70){\circle*{8}}
\put(14,40){\line(1,0){202}}
\put(190,44){\line(0,1){22}}
\put(70,70){\circle*{8}}
\put(70,40){\line(0,1){30}}
\put(4,25){$-1$}
\put(37,25){$0$}
\put(67,25){$1$}
\put(97,25){$2$}
\put(127,25){$3$}
\put(157,25){$4$}
\put(187,25){$5$}
\put(217,25){$6$}
\put(82,68){$7$}
\put(202,68){$8$}
\end{picture}
\caption{\label{d8pp}\sl Dynkin diagram of $DE_{10}\equiv D_8^{++}$
 with numbering of nodes.}
\end{figure}
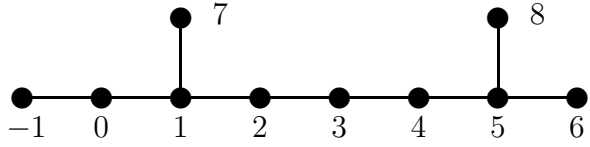

For $B_8$ the simple roots are (as labeled in fig.~\ref{b8pp})
\be\label{B8simple}
B_8:\quad\quad\quad a_1 = e_3 \quad &,& \quad 
a_2 = \frac12\big( -e_1 -e_2 - e_3 + e_4 \big) \nn\\
a_3 = e_1 \quad &,& \quad 
a_4 = \frac12\big( - 1 -e_1 - e_4 + e_5 \big) \nn\\
a_5 = 1 \quad &,& \quad 
a_6 = \frac12\big( - 1 - e_5 - e_6 - e_7 \big) \\
a_7 = \frac12\big( e_2 - e_3 +e_6 -e_7 \big)\quad &,& \quad 
a_8 = \frac14\big(  e_2 + e_4 + e_5 +e_6 +2 e_7 \big) \,.\nn
\ee
Compared to the $D_8$ simple roots of (\ref{D8simple}) the only difference is
in the simple root $a_8$, which is no longer an octavian, in agreement with
the fact that $B_8$ is not a subalgebra of $E_8$. A specific linear
combination of the $B_8$ simple 
root gives the $D_8$ expressions, describing an
embedding of $D_8$ into $B_8$. The highest root of $B_8$ is given by
\be
\theta^B &=& 2a_1 + 2 a_2 + 2 a_3 + 2 a_4 + 2 a_5 +  2a_6 +  a_7 +  2a_8 \nn\\
      &=& \frac12 \big( e_3 + e_4 + e_5 - e_7 \big) \,.
\ee
which is identical to the highest $E_8$ and $D_8$ roots.

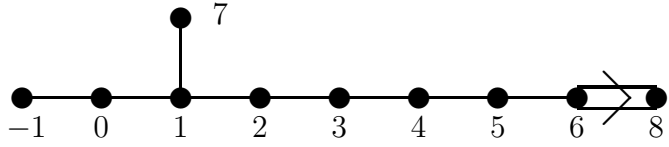
\begin{figure}
\centering
\begin{picture}(260,80)
\thicklines
\multiput(10,40)(30,0){9}{\circle*{8}}
\put(14,40){\line(1,0){202}}
\put(220,44){\line(1,0){30}}
\put(220,36){\line(1,0){30}}
\put(230,50){\line(1,-1){10}}
\put(230,30){\line(1,1){10}}
\put(70,70){\circle*{8}}
\put(70,40){\line(0,1){30}}
\put(4,25){$-1$}
\put(37,25){$0$}
\put(67,25){$1$}
\put(97,25){$2$}
\put(127,25){$3$}
\put(157,25){$4$}
\put(187,25){$5$}
\put(217,25){$6$}
\put(82,68){$7$}
\put(247,25){$8$}
\end{picture}
\caption{\label{b8pp}\sl Dynkin diagram of $BE_{10}\equiv B_8^{++}$
 with numbering of nodes.}
\end{figure}

\end{subsection}

\begin{subsection}{Weyl group $W(E_8)$}

It is known that the Weyl group of $E_8$ is of order $696\,729\,600$ and has the
structure \cite{Atlas}
\be
W(E_8) = 2\,\bd \, O_8^+(2)\,\bd \, 2 \,.
\ee
Since $W(E_8)$ naturally is a discrete subgroup of the continuous
$O(8;\reals)$ of real dimension $28$, we first record some facts about
$O(8;\reals)$. The space of unit octonions is
isomorphic to the seven sphere $S^7$ and hence of dimension $7$. At the level
of Lie algebras one has \cite{Su84,Ba02}
\be\label{so8}
\mathfrak{so}(8;\reals) \cong 
  \text{Im}(\oct) \oplus  G_2 \oplus  \text{Im}(\oct) 
\ee
as a direct sum of vector spaces but not as Lie algebras. Here,
$\text{Im}(\oct)$ denotes the purely imaginary octonions which constitute the
tangent space to $S^7$ at the octonion unit $1$. The exceptional
$14$-dimensional Lie algebra $G_2$ is isomorphic to the algebra of derivations
of the octonion algebra; the associated Lie group, which we also denote by
$G_2$, is the group of automorphism of the octonionic multiplication table. The
formula (\ref{so8}) suggests a representation of an arbitrary finite
$SO(8;\reals)$ transformation as constructed from two unit octonions
$(a_L,a_R)$ and a $G_2$ automorphism $\gamma$. The following formula
was proposed in \cite{Ra77} 
\be\label{so8act}
z \mapsto (a_L \gamma(z)) a_R
\ee
where 
\be\label{g2trm}
\gamma(z) = \bar{g}_1\bar{g}_2 (g_1 (g_2 z \bar{g}_2) \bar{g}_1) g_2g_1
\ee
is also expressed in terms of unit octonions $g_1$ and $g_2$. This formula
reproduces the correct infinitesimal $\mathfrak{so}(8,\reals)$ transformations
when 
all unit octonions are close to the identity and expressed as $a_L =  
\exp(\alpha_L)$, etc., for some $\alpha_L\in\text{Im}(\oct)$. However, it is not 
true that any finite transformation is of this form as we will show below, see
also \cite{MaDr99}. It is known that any $SO(8;\reals)$ transformation can be
written as seven times iterated (left) multiplication of $z$ by unit
octonions \cite{CoSm03}. 

Turning to the finite $E_8$ Weyl group we observe that a statement similar to
(\ref{so8}) is true for the integral octavians $\Oo$. For this we first need
some facts about the discrete automorphism group of $\Oo$. 

The group $Aut(\Oo)$ is a finite group of order $12\,096$,
usually denoted by $G_2(2)$ \cite{CoSm03}. This group is not simple
but has as simple part $U_3(3)$ of order $6\,048$. One obtains an explicit
description of a generating set of $U_3(3)$ in the following way. 
As above, the transformation
$z\mapsto a z \bar{a}$ (for $a,z\in \Oo$) gives an automorphism
iff $a^3\in \reals$, i.e. $a^3=\pm 1$. Such $a\in \Oo$ were called {\em
 Brandt transformers} in \cite{Zo35}.  Among the simple roots of $E_8$ given in
(\ref{E8simple}) this is true for the simple roots
$a_4,a_5,a_6,a_8$ which constitute a $D_4$ tree in the $E_8$
diagram. The automorphisms of $\Oo$ generated by $z\mapsto
a_iz\bar{a}_i$ (for $i=4,5,6,8$) generate in fact all of $U_3(3)$. In
order to describe the index two extension needed for the full $G_2(2)$
it is instructive to consider the following four `diagonal'
automorphisms of $z=\sum_{i=0}^7 n_i e_i$ (where $e_0=1$) 
\be
(n_0,n_1,n_2,n_3,n_4,n_5,n_6,n_7)  &\mapsto&
(n_0,n_1,n_2,n_3,n_4,n_5,n_6,n_7)\,,\nn\\ 
(n_0,n_1,n_2,n_3,n_4,n_5,n_6,n_7)  &\mapsto&
(n_0,-n_1,-n_2,n_3,n_4,-n_5,n_6,-n_7)\,,\nn\\ 
(n_0,n_1,n_2,n_3,n_4,n_5,n_6,n_7)  &\mapsto&
(n_0,n_1,-n_2,n_3,-n_4,-n_5,-n_6,n_7)\,,\nn\\ 
(n_0,n_1,n_2,n_3,n_4,n_5,n_6,n_7)  &\mapsto&
(n_0,-n_1,n_2,n_3,-n_4,n_5,-n_6,-n_7)\,,\nn 
\ee
corresponding to choosing the three lines passing through unit $e_3$
as associative triples for the Dickson doubling process.\footnote{Dickson 
 doubled pairs $(a,b)$ of quaternions manifestly have the
 automorphism $(a,b)\mapsto (a,-b)$.} This suggests that one can
obtain another automorphism by choosing a line not passing through
$e_3$, for example 
\be
(n_0,n_1,n_2,n_3,n_4,n_5,n_6,n_7)  &\mapsto&
(n_0,-n_1,-n_2,-n_3,n_4,n_5,-n_6,n_7)\,.\nn 
\ee
This automorphism is a $\ints_2$ and one can check that it combines
with $U_3(3)$ above to give all automorphisms of $\Oo$. In this way
one obtains eight diagonal automorphisms: one is the identity, while the other seven
correspond to treating any of the seven lines of the Fano plane as
the quaternion subalgebra entering the Dickson doubling process.\footnote{The
 Fano plane itself has the automorphism group 
 $PSL(2,7)$ of order $168$ and there are several copies contained in
 $Aut(\Oo)$. None is generated as a subgroup of the permutation
 group acting on the seven imaginary units.} 

Combining $Aut(\Oo)$ with two sets of octavian units (of order
$240$), one finds that there is a total of 
\be\label{ordwe8}
240\times |Aut(\Oo)|\times 240 = 696\,729\,600 = |W(E_8)|
\ee
elements, coinciding with the order of the $E_8$ Weyl group. This is the
promised discrete (and finite) version of (\ref{so8}). It is again tempting to 
associate the action of this group on $z\in\Oo$ with the expression $z\mapsto
(a_L\gamma(z))a_R$, similar to (\ref{so8act}), in particular since (\ref{so8})
is invariant under $(a_L,a_R)\leftrightarrow (-a_L,-a_R)$ thereby reducing the
tentative number of transformation described in this way to the order of the
even Weyl group $W^+(E_8)$, consistent with the fact that the odd Weyl
transformations will involve an additional conjugation of $z$.
However, the map (\ref{so8act}) has a kernel in the discrete case
in the sense that there are non-trivial triples $(a_L,\gamma,a_R)$ which act
trivially on $z\in\Oo$: Consider a $\gamma$ given by a Brandt transformer
$a\in\Oo$ with $a\neq 1$. Choosing $a_L=\bar{a}$ and $a_R=a$ leads to
$z\mapsto z$ by alternativity of $\Oo$. There are $56$ Brandt transformers
among the $240$ units of $\Oo$ and one can check that the kernel of the action
(\ref{so8act}) comes solely from situations of the type just described. Since
we do not have the correct formula for any $w\in W(E_8)$ expressed in terms of
unit octavians and a $G_2(2)$ transformation (also expressed in terms of unit
octavians) we cannot give a closed description of $W(E_8)$ using the set
(\ref{ordwe8}). This will also impede our giving a fully explicit description
of $W(E_8^{++})$. 

We close this section by explaining why the formula for finite $G_2$
automorphisms given in (\ref{g2trm}) is not correct for arbitrary unit
octonions 
$g_1$ and $g_2$. Using the Moufang identities one can show that $z\mapsto g_1
(g_2 z \bar{g}_2) \bar{g}_1$ is an automorphism iff $g_1^2 g_2^3g_1\in
\reals$, so in particular whenever $g_1$ and $g_2$ are Brandt
transformers. Conjugating with $\bar{g}_1\bar{g}_2$ to obtain (\ref{g2trm})
leads to the necessary and sufficient criterion 
\be\label{ramcrit}
\left\{\left\{\bar{g}_1, \bar{g}_2\right\}, g_2\right\} \in \reals
\ee
for (\ref{g2trm}) to be an automorphism. Here, $\left\{x,y\right\}=
xyx^{-1}y^{-1}$ is the (group) commutator which is well defined by virtue of
alternativity of $\oct$. The criterion (\ref{ramcrit}) is clearly not
satisfied for arbitrary octavians; one example is given by $g_1=a_2$ and
$g_2=a_1$ in terms of the $E_8$ simple roots. We note as a curiosity
that $W^+(E_8)$ has a subgroup of type $H_4$ which is
non-crystallographic and plays a role in the theory of quasi-crystals
\cite{MoPa93}.

\end{subsection}

\begin{subsection}{$W^+(E_{10})\equiv W^+(E_8^{++})$}

As already remarked after equation (\ref{evenmat}) the formula $S\capX
S^\dagger$ for general even Weyl transformations of the even hyperbolic Weyl
group 
acting on Hermitian $(2\times 2)$ matrices ceases to be valid in the
octonionic case. The even Weyl transformations form a discrete subgroup of
$SO(9,1;\reals)$ and similar to (\ref{so8}) it is known that at the level of
Lie algebras \cite{Su84,Ba02}
\be\label{so19}
\mathfrak{so}(1,9;\reals) \cong L'_2(\oct) \oplus \text{Im}(\oct) \oplus G_2
\;\;(\cong \mathfrak{sl}_2(\oct)\;)\,,
\ee
again as a sum vector spaces. Here, $L_2'(\oct)$ denotes the $24$-dimensional
vector space of all octonionic traceless $(2\times 2)$ matrices. Combining
this with $\text{Im}(\oct)$ one obtains the $31$-dimensional vector space of
all octonionic $(2\times 2)$ matrices with vanishing real part of the
trace. This is analogous to (\ref{pso15}) in the quaternionic case; the
additional complication of non-associativity is that one also requires an
`intertwining' $G_2$ automorphism of $\oct$. 

Similar to the suggestive, but incorrect, formula (\ref{so8act}) one could
envisage the action of $PSO(1,9;\reals)$ on $H_2(\oct)$ to be given by
\be\label{so19act}
\capX \mapsto S\gamma(\capX) S^\dagger \,.
\ee
However, this `definition' is ambiguous since it requires a prescription
for how to place the parentheses in case there are more than two
independent octonionic entries involved. Furthermore, we expect 
that in analogy with the finite $E_8$ case
discussed above that (\ref{so19act}) is not general enough to describe all
transformations. Presumably both shortcomings of (\ref{so19act}) are related
and can be resolved together. We stress nevertheless that Theorem~\ref{EvenW}
is still applicable and any even Weyl transformation can be described by an
iterated action of $(2\times 2)$ matrices.

Though we cannot offer a complete resolution to this problem we make a
comment on the issue of placing parentheses. In \cite{Ma93} the approach 
was taken that $S$ has to be such that there is no ambiguity when placing 
parentheses in the matrix expressions. This leads to very restrictive 
conditions on $S$, allowing essentially only one octonionic entry, which 
is the situation covered by Theorem~\ref{Weylprop}. However, we are here
interested in the case with more than one independent octonionic entry.
In order to make sense of (\ref{so19act}), one must presumably define
the triple product by putting parentheses {\em inside the matrix elements
in such a way that the resulting matrix is again Hermitian} (simply 
bracketing whole matrices will not do). This is for example required
to reduce (\ref{so19act}) to (\ref{so8act}) for the embedding of $W(E_8)$ 
in $W(E_{10})$ for the cases when (\ref{so8act}) is correct.

Defining the Weyl group of $E_{10}$ by iterated action of the basic
matrices given in Theorem~\ref{Weylprop} with the simple roots 
(\ref{E8simple}) (for which there arise no ambiguities due to 
non-associativity) we arrive at 
\be\label{we10}
W^+(E_8^{++}) \cong PSL_2(\Oo)\,,
\ee
where the group $PSL_2(\Oo)$ on the r.h.s. is {\em defined} by the iterated 
action. Since, unlike the $D_4$ Dynkin diagram, the $E_8$ diagram has no outer 
automorphism no additional quotients of $PSL_2(\Oo)$ are necessary for
describing the Weyl group. It is an outstanding problem to find a better
and more `intrinsic' definition of the group $PSL_2(\Oo)$, and to explore 
its implications for an associated theory of modular forms. We point out
that $PSL_2(\Oo)$ has a rich structure of subgroups, of types $\PGL$, 
$PSL_2(\cE)$, $W(D_9)$, and others, which remains to be exploited.
We hope that the information on these subgroups we have obtained in 
the preceding sections will help to find a better description of (\ref{we10}). 

\end{subsection}

\begin{subsection}{$W^+(D_8^{++})$ and $W^+(B_8^{++})$}

We close by giving the relations between the Weyl groups $W^+(D_8^{++})$ and
$W^+(B_8^{++})$, and  $W^+(D_8^{++})$ and $W^+(E_8^{++})$, respectively.

One can check that 
\be
 S_8^D=S_8^B\circ S_6^B \circ S_8^B \quad \text{and}\quad (S_8^D)^2 =\id
\ee
by using either the abstract relations or iterated matrix action of octonionic
matrices on $\capX$; the $\circ$
is meant to indicate that the product is to be understood as an iterated
action. All other generators are identical and since $S_8^B$ is acting by a
conjugation automorphism this shows that 
\begin{prop}
$W^+(B_8^{++}) \cong W^+(D_8^{++})\rtimes 2$. 
\end{prop}

Turning to $W^+(D_8^{++})$ one finds that $S_7^D$ is the only generator which
is not common to $W^+(D_8^{++})$ and $W^+(E_8^{++})$ since all other simple
roots and the highest roots are identical. A calculation reveals that $S_7^D$
can be expressed as the result of a conjugation action of the $W^+(E_8^{++})$
generators as 
\be\label{d8e8}
S_7^D&=&S_7 S_6 S_5 S_8 S_4 S_3 S_5 S_4 S_6 S_7 S_5 S_8 S_6 S_5 S_4 S_3
S_2\nn\\  
  &&\circ S_3 S_4  S_5 S_6 S_8 S_5 S_7 S_6 S_4 S_5 S_3 S_4 S_8 S_5 S_6 S_7\,,
\ee
where we omitted a superscript on $S^E_i$ and the $\circ$ symbols on the
r.h.s. in order not to clutter notation. The Weyl transformation (\ref{d8e8})
maps the simple root $\a_7^E$ to $\a_7^D$.
This shows that $W^+(D_8^{++})$ is a
subgroup of $W^+(E_{10})$ but one can check that it is not normal in
$W^+(E_{10})$ since conjugation by $S_7^E$ does not preserve the subgroup.
For the finite Weyl groups $W(D_8)$ has index $135$ in $W(E_8)$.
By contrast, for the hyperbolic Weyl groups one has the following result from~\cite{Johnson99,Tumarkin04}.
\footnote{We thank N.~Johnson and P.~Tumarkin for bringing to our attention a mistake in the published 
version of the following proposition.}

\begin{prop}
$W^+(D_8^{++})$ is an index $527$ subgroup of
$W^+(E_8^{++})\cong PSL_2(\Oo)$.
\end{prop}

Without a direct definition of $PSL_2(\Oo)$ we cannot give a more
detailed octonionic description of $W^+(D_8^{++})$ that would be
analogous to (\ref{gq}). There is no direct relation between 
$W^+(B_8^{++})$ and $W^+(E_8^{++})$.

\end{subsection}

\end{section}


\medskip
\appendix

\begin{section}{Examples with twisted affine algebras}
\label{twapp}

In this appendix we show that our techniques are also suited for treating
hyperbolic algebras which do not arise as over-extensions but whose Dynkin
diagrams instead involve subdiagrams of twisted affine type. For simplicity we
will exemplify this for two cases constructed over the complex numbers with
twisted affine algebras $D_2^{(2)}$ and $D_4^{(3)}$ in the terminology of
\cite{Ka90}, both of which have rank three. The hyperbolic node is attached to
the `twisted affine' node with the single node for the cases we consider and
we therefore denote the associated hyperbolic algebras by $D_2^{(2)+}$ and
$D_4^{(3)+}$. However, we anticipate that there are also other, more general
hyperbolic cases which fit into the picture we develop in this paper.

\begin{subsection}{$\Kk=\cx$, type $D_2^{(2)+}$}

The hyperbolic Kac--Moody algebra $D_2^{(2)+}$ can be constructed from the
root system of type $B_2$. Although the root system is identical to that of
$C_2$ already depicted in fig.~\ref{c2roots} we draw it again in
fig.~\ref{b2roots} where we highlight the short highest root which enters the
twisted affine extension $D_2^{(2)}$. 

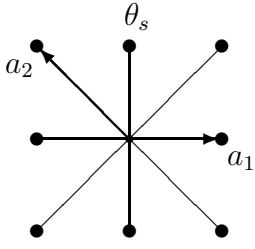
\begin{figure}
\centering
\begin{picture}(90,90)
\put(10,10){\line(1,1){70}}
\put(10,45){\line(1,0){70}}
\put(10,80){\line(1,-1){70}}
\put(45,10){\line(0,1){70}}
\put(45,45){\circle*{3}}
\put(10,10){\circle*{5}}
\put(10,45){\circle*{5}}
\put(45,10){\circle*{5}}
\put(80,10){\circle*{5}}
\put(45,80){\circle*{5}}
\put(43,88){$\theta_s$}
\put(80,45){\circle*{5}}
\put(82,35){$a_1$}
\put(80,80){\circle*{5}}
\put(10,80){\circle*{5}}
\put(-2,70){$a_2$}
\thicklines
\put(45,45){\vector(1,0){34}}
\put(45,45){\vector(-1,1){34}}
\end{picture}
\caption{\label{b2roots}\sl The root system of type $B_2$, with simple roots
  labeled and indicated by arrows. The lattice they generate is the ring of
  Gaussian integers. There are both long and short roots and we indicate the
  highest short root $\theta_s$.} 
\end{figure}

The root system is not simply laced, having simple roots whose squared lengths
are in the ratio $2$ to $1$. We choose the following units and simple roots
for $B_2$ in the complex plane 
\be\label{typeB2}
\ve_1 = a_1 =  1 \;\; , \;  \ve_2 = \frac{a_2}{\sqrt{2}} = \frac{-1 +
  i}{\sqrt{2}} 
  \;\; , 
  \; \theta_s = i\,.
\ee
These formulas are identical to those of (\ref{typeC2}) except that the whole
lattice has been rescaled by a factor $\sqrt{2}$. 

The hyperbolic simple roots of $D_2^{(2)+}$ now take the form
\begin{align}
\a_{-1}&=\lp\begin{array}{cc}1&0\\0&-1\end{array}\rp\,,&
\a_0&=\lp\begin{array}{cc}-1&-\theta_s\\-\bar\theta_s&0\end{array}\rp = 
\lp\begin{array}{cc}-1&-i\\i&0\end{array}\rp \,,&\nn\\
\a_{1}&=\lp\begin{array}{cc}0&1\\1&0\end{array}\rp\,,&
\a_{2}&=\lp\begin{array}{cc}0&-1+i \\-1-i &0\end{array}\rp\,.&
\end{align}
The Dynkin diagram of the hyperbolic algebra $D_2^{(2)+}$ is shown in
fig.~\ref{b2pp}. The non-zero inner products between the simple roots are 
\be\label{b2cn1}
(\a_{-1},\a_{-1}) = (\a_0,\a_0) = (\a_1,\a_1) = 2\;\; , \; (\a_2,\a_2) =
4\;,\nn 
\\\label{b2cn2}
(\a_{-1},\a_0) = -1\;\; , \; (\a_0,\a_2) = -2\;\; , \; (\a_1,\a_2) = -2
\ee
and these give rise to the $D_2^{(2)+}$ diagram.

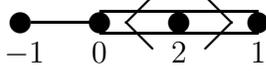
\begin{figure}
\centering
\begin{picture}(110,20)
\thicklines
\multiput(10,10)(30,0){4}{\circle*{8}}
\put(14,10){\line(1,0){22}}
\put(40,14){\line(1,0){30}}
\put(40,6){\line(1,0){30}}
\put(70,14){\line(1,0){30}}
\put(70,6){\line(1,0){30}}
\put(80,20){\line(1,-1){10}}
\put(80,0){\line(1,1){10}}
\put(60,20){\line(-1,-1){10}}
\put(60,0){\line(-1,1){10}}
\put(4,-5){$-1$}
\put(37,-5){$0$}
\put(67,-5){$2$}
\put(97,-5){$1$}
\end{picture}
\caption{\label{b2pp}\sl $D_2^{(2)+}$ Dynkin diagram with numbering of
 nodes.}
\end{figure}

Turning to the Weyl group of $D_2^{(2)^+}$,  we realize that the Weyl group is
isomorphic to that of $C_2^{++}$ since the Coxeter group is insensitive to the
direction of the arrows. 
\begin{prop}
$W^+(D_2^{(2)+}) \cong W^+(B_2^{++}) \cong PSL_2(\cG) \rtimes 2$.
\end{prop}

\end{subsection}

\begin{subsection}{$\Kk=\cx$, type $D_4^{(3)+}$}

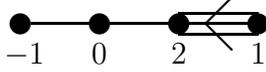
\begin{figure}
\centering
\begin{picture}(110,20)
\thicklines
\multiput(10,10)(30,0){4}{\circle*{8}}
\put(14,10){\line(1,0){22}}
\put(40,10){\line(1,0){30}}
\put(70,14){\line(1,0){30}}
\put(70,10){\line(1,0){30}}
\put(70,6){\line(1,0){30}}
\put(90,20){\line(-1,-1){10}}
\put(90,0){\line(-1,1){10}}
\put(4,-5){$-1$}
\put(37,-5){$0$}
\put(67,-5){$2$}
\put(97,-5){$1$}
\end{picture}
\caption{\label{d43p}\sl $D_4^{(3)+}$ Dynkin diagram with numbering of nodes.}
\end{figure}

The final hyperbolic algebra we consider is the hyperbolic extension of the
twisted affine 
algebra of type $D_4^{(3)}$ (the only one twisted with an order three
automorphism). We denote the hyperbolic extension by $D_4^{(3)+}$  
and its diagram is shown in fig.~\ref{d43p}. Its construction employs the
$G_2$ root system, shown again in fig.~\ref{d43rts}, and we make the following
choices of units and simple roots in the complex plane: 
\be\label{typeDualG2}
\ve_1 = \frac{a_1}{\sqrt{3}} =  1  \;\; , \; \ve_2  = a_2=
\frac{-\sqrt{3}+i}{2}   
  \;\; , 
  \; \theta_s = i \,.
\ee
Again, these formulas are those of (\ref{typeG2}) except for a rescaling by a
factor $\sqrt{3}$ and the use of the highest short root $\theta_s$ instead of
the highest (long) root $\theta$. 

\begin{figure}
\centering
\begin{picture}(100,140)
\put(115,48){$a_1$}
\put(-14,75){$a_2$} 
\put(37,110){$\theta_s$}
\put(40,60){\line(0,1){40}} 
\put(40,60){\line(0,-1){40}} 
\put(40,60){\line(2,1){40}} 
\put(40,60){\line(-2,1){40}} 
\put(40,60){\line(2,-1){40}} 
\put(40,60){\line(-2,-1){40}} 
\put(80,40){\circle*{5}} 
\put(80,80){\circle*{5}} 
\put(40,20){\circle*{5}} 
\put(40,100){\circle*{5}} 
\put(0,40){\circle*{5}} 
\put(0,80){\circle*{5}} 
\put(40,60){\line(1,0){80}}
\put(40,60){\line(-1,0){80}}
\put(40,60){\line(2,3){40}}
\put(40,60){\line(-2,3){40}}
\put(40,60){\line(2,-3){40}}
\put(40,60){\line(-2,-3){40}}
\put(-40,60){\circle*{5}} 
\put(120,60){\circle*{5}} 
\put(0,120){\circle*{5}} 
\put(0,0){\circle*{5}} 
\put(80,120){\circle*{5}} 
\put(80,0){\circle*{5}} 
\thicklines
\put(40,60){\vector(1,0){80}}
\put(40,60){\vector(-2,1){40}}
\end{picture}
\caption{\label{d43rts}\sl The root system of type $G_2$ with highest  short
 root indicated.}  
\end{figure}
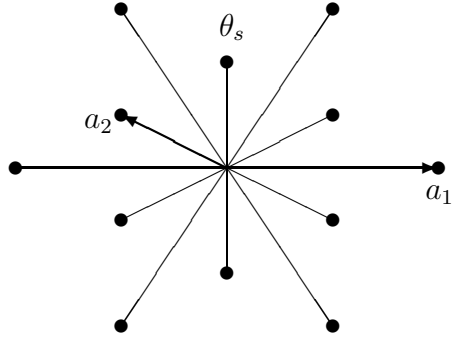

The associated hyperbolic simple roots of $D_4^{(3)+}$ are then
\begin{align}
\a_{-1}&=\lp\begin{array}{cc}1&0\\0&-1\end{array}\rp\,,&
\a_0&=\lp\begin{array}{cc}-1&-\theta_s\\-\bar\theta_s&0\end{array}\rp
  =\lp\begin{array}{cc}-1&-i\\i&0
\end{array}\rp\,,&\nn\\
\a_{1}&=\lp\begin{array}{cc}0&\sqrt{3}\\\sqrt{3}&0\end{array}\rp\,,&
\a_{2}&=\lp\begin{array}{cc}0&\frac{-\sqrt{3}+i}{2}\\
\frac{-\sqrt{3}-i}{2}&0\end{array}\rp\,.& 
\end{align}
The non-zero inner products between these simple roots are
\be
(\a_{-1},\a_{-1}) = (\a_0,\a_0) = (\a_2,\a_2) = 2\;\; , \; (\a_1,\a_1) =
6\,,\nn\\ 
(\a_{-1},\a_0) = -1\;\; , \; (\a_0,\a_2) = -1\;\; , \; (\a_1,\a_2) = -3\,.
\ee

The Weyl group is again easy to determine in this case since it is
isomorphic to that of the standard $G_2^{++}$ case. 
\begin{prop}
$W^+(D_4^{(3)+}) \cong W^+(G_2^{++})\cong PSL_2(\cE) \rtimes 2$.
\end{prop}

\end{subsection}

\end{section}

\end{document}